\documentclass[12pt,draftcls,onecolumn]{IEEEtran}

\IEEEoverridecommandlockouts                              % This command is only
\def\1{{\bf 1}}
\def\R{\mathbb{R}}

\usepackage{color}
\usepackage{algorithmic}
\usepackage{ amssymb }
\usepackage{epsfig}
\def\ao#1{{#1}}
\def\aoj#1{{#1}}
\def\aoc#1{{#1}}
\def\aor#1{{#1}}
\title{\LARGE \bf
Nonuniform coverage control on the line
}

\author{Naomi Ehrich Leonard and Alex Olshevsky % <-this % stops a space
\thanks{Naomi Leonard is with the Department of Mechanical and Aerospace Engineering, Princeton University, Princeton, NJ
08544, USA, {\tt \small naomi@princeton,edu}}
\thanks{Alex Olshevsky is with the Department of Industrial and Enterprise Systems Engineering, University of Illinois at
Urbana-Champaign, Urbana, IL, 61801, USA, {\tt \small aolshev2@illinois.edu}}%
\thanks{This research was supported in part by AFOSR grant FA9550-07-1-0-0528 and ONR grant N00014-09-1-1074.}
}

\begin{document}

\maketitle
\thispagestyle{empty}
\pagestyle{empty}

%%%%%%%%%%%%%%%%%%%%%%%%%%%%%%%%%%%%%%%%%%%%%%%%%%%%%%%%%%%%%%%%%%%%%%%%%%%%%%%%
\begin{abstract} This paper investigates control laws allowing mobile, autonomous agents
to optimally position themselves on the line for distributed sensing in
a nonuniform field. We show that a simple static control law, based
only on local measurements of the field by each agent, drives the
agents close to the optimal positions after \aor{the agents execute in parallel} a
number of sensing/movement/\aor{computation} rounds that is essentially quadratic in the
number of agents. Further, we exhibit a dynamic control law which,
under slightly stronger assumptions on the capabilities and knowledge
of each agent, drives the agents close to the optimal positions after \aor{the agents execute in parallel} a number
of sensing/communication/\aor{computation}/movement rounds that is essentially linear in
the number of agents. Crucially, both algorithms are
fully distributed and robust to unpredictable loss and addition of
agents.
\end{abstract}

%%%%%%%%%%%%%%%%%%%%%%%%%%%%%%%%%%%%%%%%%%%%%%%%%%%%%%%%%%%%%%%%%%%%%%%%%%%%%%%%
\section{INTRODUCTION}

Widespread deployment of networks of sensors and autonomous vehicles is
expected to revolutionize our ability to monitor and control physical
environments from remote locations. However, for such networks to achieve their
full range of applicability, they must be capable of operating in uncertain 
and unstructured environments without centralized supervision. Realizing the
full potential of such systems will require
the development of protocols that are fully autonomous, distributed,
and adaptive in the face of  changing environments.

An important problem in this context is the coverage problem:  a network of mobile sensors 
should distribute itself over a region with a measurable field  so that the likelihood of detecting an event of interest is maximized. 
If the probability distribution of the event is uniform over the area, then the optimal solution
will involve a uniform spatial distribution  of the agents. On the other hand if this probability distribution
is not uniform, then the sensors should be more densely positioned in the subregions that
have higher event probability. 

\aoj{There has been considerable recent interest in the coverage problem, in part prompted by applications 
to the design of robotic networks for search and rescue, environmental monitoring, automatic surveillance, cooperative estimation, and other tasks (see survey below). However, it is still unknown how to design decentralized control  laws for achieving optimal coverage which scale optimally or nearly-optimally with the number of agents in the network.}  \aoj{In this paper, we do this for the nonuniform coverage problem in the case the agents are arranged on a line. Moreover, we establish that relatively minor increases in the capabilities of each agent can lead to an order of magnitude improvement in the convergence rate of the system in achieving coverage. }

\subsection{Related work}

There is a considerable literature on coverage algorithms for groups of dynamic agents; we refer the reader to \cite{hms02, zk03, ps04, cmkb04, cb05, cmb05, lc05, AKPJK06, AJKJ07, hs07, hss07b, gg07, hss07, AKPJK07, ghhf08, wtk08, pkmp08, sbsr08, af08, ch08a, ch08b, zc08,  bmh09, zc09, ll09, ssr09, pt09, bsmsr10}
and the references therein. \aoj{Here, we only survey the papers most closely related to our work.}

In \cite{cb05}, uniform coverage algorithms are derived using Voronoi cells and gradient laws for
distributed dynamical systems. Uniform constrained coverage control is addressed in
\cite{ps04} where the constraint is a minimum limit on node degree. Virtual potentials enable repulsion between agents to maximize
coverage and attraction between agents to enforce the constraint. In \cite{lc05}, gradient
control laws are proposed to move sensors to a configuration that maximizes expected
event detection frequency. Local rules are enforced by defining a sensing radius for
each agent, which also makes computations simpler. The approach is demonstrated
for a nonuniform but symmetric density field with and without communication constraints. Further results for distributed coverage control are presented in \cite{cmkb04} for a coverage metric defined in terms of the Euclidean metric with a weighting factor that
allows for nonuniformity. As in \cite{cmkb04}, the methodology makes use of Voronoi cells and
Lloyd descent algorithms. \ao{The papers \cite{ch08a,ch08b} identified 
a class of non-convex regions for which the coverage problem may be solved by reduction to the convex case through a well-chosen transformation of the
region. The papers \cite{AKPJK06,AJKJ07} explored an optimization-based approach to some complex variations of the covering problem.  } 

The paper \cite{ll09} considered the general nonuniform coverage problem with a
non-Euclidean distance, and it proposed and proved the correctness of a coverage control law in the plane.
However, the control law of \cite{ll09} is only partially distributed, in that it relies on a ``cartogram
computation'' step which requires some global knowledge of the domain.

\aor{ \subsection{The coverage problem on the line} We will consider a particular case of the nonuniform coverage problem when the agents are arranged on a 
line. There are two main reasons for our interest in this question. First, while this is arguably the simplest model within which one can study the coverage problem, even in this setting many of the most fundamental questions are open: for example, it is unclear how to design control laws for which the time taken or energy expended scales optimally with the number of agents. }

\aor{Secondly, the line coverage problem is motivated by several natural applications. One is the problem of border patrol, in which a collection of mobile sensors must adaptively spread themselves out over a one-dimensional curve; this is easily seen to be equivalent to the problem we consider here after a re-parametrization.  Boundary tracking is a problem analogous to border patrol; in this case mobile sensors should adapt their positions to monitor the boundary of a crowd of people, a toxic spill, a forest fire, etc.  Another application is  adaptive sampling in two and three-dimensional environments where the nonuniformity in the sampled field is dominant in one dimension.   For instance, in the air or in the water, nonuniform coverage in the vertical is a typical concern:  vertical profilers or autonomous vehicles should distribute themselves over depth or altitude, while vertical lines of these sensors can be deployed in the horizontal plane in a static array or adapted using a distributed uniform coverage algorithm.   An important example in the vertical water column (oceans, lakes, reservoirs, rivers) is mapping the  nonuniform vertical temperature field and locating the thermocline, a thin  layer that separates the warm surface layer from the cold deep layer; the depth of the thermocline is relevant in many applications, for instance, it is the depth at which the speed of sound is maximal.  
}     

\subsection{Our contributions} Our work  builds on the results of \cite{ll09} to design a control law for the nonuniform
coverage problem in the one-dimensional case when the agents are positioned on the 
line. We develop fully distributed coverage control laws for a nonuniform field in this setting, and moreover, 
we prove 
quantitative convergence bounds on the performance of these algorithms. Interestingly,
we find that relatively modest increases in the capabilities and knowledge of each agent
can translate into considerable improvements in the global performance. \ao{These improvements are 
obtained by implementing more sophisticated variants of distributed algorithms; in particular, we 
make significant use of the technique of lifting Markov chains, introduced in \cite{dhn00}, and subsequently
used to accelerate distributed computation in \cite{jss10,s09}).}

We begin with an introduction 
to the nonuniform coverage problem in Section II. In Section III, we present our first fully distributed
control law for the coverage problem. The execution of this control law only requires the agents 
to be able to measure distances to their neighbors and to measure the field around their 
location. The main result of this section is Theorem 1, which demonstrates
the correctness of the algorithm and gives a quantitative bound on its performance. We show that 
it takes $n$ agents essentially on the order of $n^2$ \aor{parallel update rounds consisting of sensing, computation,
and movement} to come close to the optimal configuration regardless of the initial conditions.

In Section IV, we present another fully distributed control law for coverage. The execution of this control
law requires more capabilities on the part of the agents: they store several numbers in memory, communicate
these numbers to their neighbors at every round, and moreover, they know approximately (within
a constant factor) how many agents there are in total. Subject to these assumptions, we derive
a considerable \aor{improvement in the number of update rounds}  over the simple static control law of Section III. The main result of this section is 
Theorem 6, which demonstrates the correctness of the algorithm and gives a quantitative bound on its performance.   We show that it takes $n$
agents essentially on the order of $n$ \aor{parallel update rounds consisting of sensing, computation, communication, and
movement} to come close to the optimal configuration regardless of initial conditions. This is an order of magnitude improvement over the
control law of Section III. 

%%%%%%%%%%%%%%%%%%%%%%%%%%%%%%%%%%%%%%%%%%%%%%%%%%%%%%%%%%%%%%%%%%%%%%%%%%%%%%%%
\section{NONUNIFORM COVERAGE}

We introduce the nonuniform coverage problem in this section; our exposition closely follows the
expositions of \cite{ll09,cb05}. We consider $n$ mobile agents initially situated at arbitrary positions $x_1(0), x_2(0), \ldots, x_n(0)$ which, for simplicity, we henceforth assume to be 
located in the interval $[0,1]$. There is a positive, piecewise-continuous function $\rho: [0,1] \rightarrow (0, \infty)$, which measures the density of information
or resource at each point. \aor{We assume that $\rho$ is bounded both above and below, i.e., there exist positive numbers $\rho_{\rm max}$ and $\rho_{\rm min}$ such that for all $z \in [0,1]$, we have
$ \rho_{\min} \leq \rho(z) \leq \rho_{\rm max}.$} The goal is to bring the agents from their initial configuration to a static configuration that allows them to optimally sense 
in the density field $\rho$. Intuitively, we would like more agents to be positioned in areas where $\rho$ is high, and fewer agents positioned in areas where $\rho$ is low. 

More formally, for $a,b \in [0,1]$ we define the metric 
\[ d_{\rho}(a,b) = \int_{\min(a,b)}^{\max(a,b)} \rho(z) dz.\] It is easy to see this defines a valid metric between points in $[0,1]$. Relative to the ordinary distance 
$|a-b|$, this metric expands regions where $\rho$ is large and shrinks regions where $\rho$ is small.

Following \cite{ll09}, we define the coverage of a set of points $x_1,\ldots,x_n$ relative to the density field $\rho$ as 
\[ \Phi(x_1, \ldots, x_n, \rho) = \max_{y \in [0,1]} \min_{i=1,\ldots,n} d_{\rho}(y,x_i). \] Given the positions $x_1, \ldots, x_n$ of the
agents and the density field $\rho$, computing $\Phi$ requires computing the distance $d_{\rho}$ from any point in $[0,1]$ to the closest $x_i$. The
coverage metric $\Phi$ is then the largest of these distances. A smaller $\Phi$ implies the vehicles achieve better coverage of the domain
$[0,1]$.  We use $\Phi^*$ to denote the best (smallest) possible coverage 
\[ \Phi^* = \inf_{(x_1, \ldots, x_n) \in [0,1]^n} \Phi(x_1,\ldots,x_n, \rho).\]  \aor{It is not hard to see (and we will formally show later) that this infimum is reached at a point 
where no two of the $x_i$ are the same.} In this paper, we are concerned with designing control laws which drive agents towards positions with coverage $\Phi^*$. 
%We assume that agents 
%do not know their own positions, nor do they know the positions of their neighbors. They can, however, measure relative distances between themselves
%and their neighbors; moreover the leftmost agent can measure the distance between itself and the ``border point'' $0$, and similarly the rightmost %agent
%can measure the relative distance to the point $1$. In addition, the agents are also capable of taking measurements of the field $\rho$ in a %neighborhood of their location. Finally, in Section IV, we will make some stronger assumptions on the capabilities of each agent, such as the ability %to store numbers in memory, the ability to communicate with left and right neighbors, the ability to sense when your left and right neighbors
%haved moved, and, finally, approximate knowledge of the total number of agents (up to a constant factor). 

As pointed out in \cite{ll09}, \aor{one may alternatively interpret the coverage problem in terms of optimal positioning} with a nonuniform distance; \aor{this} is closely related to information gathering and sensor array optimization problems. A typical problem
is to minimize shortest response time from a collection of vehicles to any point in a terrain of varying roughness. In that case,  the non-Euclidean distance $d_{\rho}$ appears
because rougher bits of terrain take longer to traverse. Another such problem is the detection of acoustic signals; the objective is to place sensors so they can
detect a source anywhere in an inhomogeneous medium. In that case, the non-Euclidean distance $d_{\rho}$ appears as a result of the spatially varying refractive index of the environment.

\section{A STATIC COVERAGE CONTROL LAW}

We now describe and analyze a simple distributed control law that drives the vehicles towards optimal coverage. 
First, we need to define the notion
of a $\rho$-weighted median between points.

\smallskip

\noindent {\bf Definition:} The $\alpha$-median $m_{\rho}^{\alpha}(a,b)$ is defined as the point 
$c \in (a,b)$ which satisfies \[ \int_a^c \rho(z) dz = \alpha \int_c^b \rho(z) dz.\] Due to the strict positivity 
of $\rho$, it is easy to see that a unique such point exists for any $\alpha \geq 0$. 

\smallskip

We can now state the coverage control law. We assume for convenience that agents are labeled $1,\ldots,n$ from left to 
right. This makes it easier to state what follows; however, the actual implementation of the algorithm does not require
the use of these labels. 

\bigskip

\noindent {\bf A static coverage control law:} the agents iterate as
\begin{eqnarray*} x_1(t+1) & = & m_{\rho}^{1/2}(0,x_2(t))  \\
x_i(t+1) & = &  m_{\rho}^1 (x_{i-1}(t), x_{i+1}(t)), ~~~ i = 2, \ldots, n-1  \\
x_n(t+1) & = & m_{\rho}^2(x_{n-1}(t),1)  
\end{eqnarray*}

\bigskip

We first briefly  outline how this scheme may be implemented without knowledge of the 
labels $1,\ldots,n$ by the nodes. A node $i$ will initially
check whether it has a left neighbor and a right neighbor, or whether it is a ``border agent'' with a single neighbor. Suppose it
has two neighbors. Then, it will measure the distance $d^{\rm L}$ to its left neighbor and $d^{\rm R}$ to its right neighbor, and denoting
its position (which it does not know) by $x_i$, it will measure $\rho$ in the interval $[x_i-d^{\rm L}, x_i+d^{\rm R}]$. This gives it enough information to compute the $1$-median of the positions of its neighbors, and it moves to this location. ``Border agents'' can similarly implement this 
control law without knowledge of their labels.

Next, we remark that this scheme may be interpreted as a distributed implementation of the 
cartogram approach introduced in \cite{ll09} specialized to the line. Intuitively, each of the middle nodes 
$2, \ldots, n-1$ seeks to position itself ``in the middle'' of its neighboring agents while 
stretching areas with high $\rho$ and shrinking areas with low $\rho$; this is precisely the distributed
computation of the cartograms used in \cite{ll09}. 

Our goal in this section is to prove that this iteration solves the coverage control problem and to provide
quantitative bounds on its 
performance. The main result of this section is the following theorem. 

\bigskip

\noindent {\bf Theorem 1:} Each $x_i(t)$ has a limit, and the limiting set of positions have coverage $\Phi^*$. Moreover, after
$O(n^2 \log (\frac{n}{\epsilon} \frac{\rho_{\rm max}}{\rho_{\rm min}}))$
%\[ 4 (n+1)^2 \log \left(  \frac{n+1}{\epsilon} \sqrt{ 2 \frac{\rho_{\rm max}}{\rho_{\rm min}}} \right) \] 
rounds, each agent is within $\epsilon$ of its final 
limit. 

\bigskip

%As far as we are aware, this is the first quantative convergence time bound for the coverage problem. 

We next turn to the proof of this theorem. We first write down the optimality conditions for achieving $\Phi^*$. 

\smallskip

\noindent  {\bf Lemma 2.} The equations \begin{small} \begin{equation} \label{optcond} 2 d_{\rho}(0,x_1) = d_{\rho}(x_1,x_2) = \cdots d_{\rho}(x_{n-1},x_n) = 2 d_{\rho}(x_n,1)\end{equation}\end{small}have a unique solution which achieves coverage $\Phi^*$.  Moreover, \[ \Phi^* = \frac{1}{2n} d_{\rho}(0,1). \] 

\noindent {\bf Proof. } It is easy to see that Eq. (\ref{optcond}) has a solution. Indeed,
suppose we position agent $1$ at the location $x_1 \in (0,1)$ satisfying $d_{\rho}(0, x_1) = d_{\rho}(0,1)/(2n)$. Such
a location clearly exists. Next, for every $i>1$, we position the $i$'th agent at
the location $x_i$ which satisfies $d_{\rho}(x_{i-1},x_i) = d_{\rho}(0,1)/n$. Again, such a point
clearly must always exist. The resulting locations satisfy Eq. (\ref{optcond}).

Moreover, this solution is unique because together with the obvious identity \[ d_{\rho}(0,1) = d_{\rho}(0,x_1) + \sum_{i=1}^{n-1} d_{\rho}(x_i, x_{i+1}) + d_{\rho}(x_n,1), \] Eq. (\ref{optcond}) implies \[ d_{\rho}(0,1)  =  d_{\rho}(0,x_1) ( 1 + 2(n-1)+1) \] which
immediately fixes the value of $x_1$, and hence all the values $x_i$. 

Under this configuration the largest distance $d_{\rho}$ 
between any point in $[0,1]$ and its closest agent clearly equals $d_{\rho}(0,1)/(2n)$. We next argue that 
$\Phi^*$ is at least this big. 
%Next, lets argue that $\Phi^*$ cannot be strictly smaller than $d_{\rho}(0,1)/(2n)$. %If it was, then there would
%be a configuration $x_1' < x_2' < \cdots x_n'$ such that $\Phi(x_1', \ldots, x_n', \rho) < d_{\rho}(0,1)/(2n)$.
Consider any set of positions $x_1 < x_2 < \cdots < x_n$. Since
\begin{small} 
\begin{eqnarray} d_{\rho}(0,1) & = &  \left[ d_{\rho}(x_1,x_2)\right] + \left[ d_{\rho}(x_2,x_3) \right] +  \cdots + \left[ d_{\rho}(x_{n-1},x_n) \right] \nonumber + \left[ d_{\rho}(0,x_1) + d_{\rho}(x_n,1) \right] \label{sumeq} 
\end{eqnarray}\end{small}and there are $n$ terms in brackets, we can conclude that at least one of the bracketed terms is at least $d_{\rho}(0,1)/n$.
 
If the only such term is the last term, then one of $d_{\rho}(0,x_1)$, $d_{\rho}(x_n,1)$ has value at least $d_{\rho}(0,1)/(2n)$. 
But this implies that the distance $d_{\rho}$ from either $0$ or $1$ to the closest agent is at least $d_{\rho}(0,1)/(2n)$, which proves that in this
case $\Phi^*$
is at least that big. 

If, on the other hand, it is some term $d_{\rho}(x_i, x_{i+1})$ which is at least $d_{\rho}(0,1)/n$, then the distance $d_{\rho}$ from 
the median $m_{\rho}^1(x_i,x_{i+1})$ to the closest agent is at least $d_{\rho}(0,1)/(2n)$, which again implies $\Phi^*$ is
at least this much. {\bf q.e.d.}

\smallskip

%The proof of Lemma 1 also proves the following. 

We next introduce a change of variables which makes our static control law easier to analyze. We define \[ F(x) = \int_0^x \rho(z) dz,\]
and note that $F(1) = d_{\rho}(0,1)$. Moreover, for any two points $a<b$ in $[0,1]$,
 \[ d_{\rho}(a,b) = F(b)-F(a),\]  and
\begin{equation} \label{F-median} F(m_{\rho}^{\alpha}(a,b)) = \frac{F(a)+\alpha F(b)}{1+\alpha}.\end{equation} %We can restate 
%Lemma 2 in terms of this mapping as follows. 

%\smallskip

%\noindent {\bf Corollary 3:} The set of points $x_1,\ldots,x_n$ achieves coverage $\Phi^*$ if and only if
%\begin{eqnarray*} F(x_1) & = & \frac{F(1)}{2n} \\
%F(x_i) & = & \frac{F(1)}{2n} + (i-1) \frac{F(1)}{n},  ~~~ i = 2,\ldots,n
%\end{eqnarray*} 
%\smallskip

\smallskip

\aor{\noindent {\bf Remark:} We note that under the execution of the static control law, we have that for all $t$, \[ x_1(t) \leq x_2(t) \leq \cdots \leq x_n(t) \] That is, the execution of the static control law preserves the 
ordering of the agents. To show this it suffices to demonstrate that for all $i,t$,
\[ F(x_i(t)) \leq F(x_{i+1}(t)) \] which follows by induction from Eq. (\ref{F-median}) and the algorithm defintion. 
%The same argument shows that \[ x_i(t + 1 ) \in [x_{i-1}(t), x_{i+1}(t)] \] and 
%similarly \[ x_1(t+1) \in [0, x_2(t)], ~~ x_{n}(t+1) \in [0, x_n(t)] \]
}

\aoj{The following simple corollary of Lemma 2 describes how the nonuniform coverage problem may be made uniform by a transformation and connects our approach to the ``nonuniform cartogram'' approach of \cite{ll09}.}

\bigskip

\aoj{\noindent {\bf Corollary 3:} The positions $x_1, \ldots, x_n$ achieve optimal coverage $\Phi^*$  
if and only if $F(x_1), \ldots, F(x_n)$ achieve optimal coverage in $[0,F(1)]$ with density equal to $1$. }

\bigskip

The next lemma restates our coverage control law in a particularly convenient 
form. 

\noindent {\bf Lemma 4:} Assuming there are at least two agents, let us define 
\begin{eqnarray*} d_0(t) & = & 2 F(x_1(t)) \\
d_i(t) & = & F(x_{i+1}(t)) - F(x_i(t)), ~~~i = 2, \ldots, n-1 \\
d_n(t) & = & 2(F(1) - F(x_n(t)))
\end{eqnarray*} and let $d(t)$ be the vector in $\R^{n+1}$ which stacks the 
variables $d_i$. Finally, we define \[ U_3 = \left( \begin{array}{ccc} 
-4 & 4 & 0 \\
2 & -4 & 2 \\
0 & 4 & -4 
\end{array} \right),\] \[ U_4 = \left( \begin{array}{cccc} -4 & 4 & 0 & 0 \\
2 & -5 & 3 & 0 \\
0 & 3 & -5 & 2 \\
0 & 0 & 4 & -4
\end{array} \right) \] and for $n \geq 5$, \[ U_{n} = \left( \begin{array}{ccccccc}
-4 & 4  &      &      &      &    &  \\
 2 & -5 & 3    &      &      &    &  \\
   & 3  & -6   & 3    &      &    &  \\
   &    &\ddots&\ddots&\ddots&    &  \\
   &    &      &3     &-6    & 3  &  \\
   &    &      &      & 3    & -5 & 2\\
   &    &      &      &      & 4  & -4
 \end{array} \right). \] Then $d(t)$ follows the dynamics
 \[ d(t+1) = (I + \frac{1}{6} U_{n+1}) d(t). \]
 
 \noindent {\bf Proof:} We  give the calculation in the case of $n \geq 4$; the proof of the 
 cases $n=2,3$ are similar. Define $y_i(t) = F(x_i(t))$.  From Eq. (\ref{F-median}), we can
 conclude that 
 the variables $y_i(t)$ evolve as \begin{eqnarray*} y_1(t+1) & = & \frac{y_2(t)}{3}  \\
y_i(t+1) & = & \frac{y_{i-1}(t) + y_{i+1}(t)}{2}, ~~~~ i = 2, \ldots, n-1 \\
y_n(t+1) & = & \frac{2F(1)+y_{n-1}(t)}{3}.
\end{eqnarray*} We can rewrite the variables $d_i(t)$ in terms of the variables
$y_i(t)$ as 
\begin{eqnarray} d_0(t) & = & 2 y_1(t) \nonumber \\
d_i(t) & = & y_{i+1}(t) - y_i(t), ~~~~i = 2, \ldots, n-1 \nonumber \\ 
d_n(t) & = & 2(F(1) - y_n(t)) \label{dintermsofy}
\end{eqnarray} Observe that
\begin{eqnarray} y_1(t+1)  & = & y_1(t) +  \frac{d_1(t)-d_0(t)}{3} \nonumber \\
y_i(t+1) & = & y_i(t) + \frac{d_i(t) - d_{i-1}(t)}{2} ~~~~ i = 2, \ldots, n-1 \nonumber \\
y_n(t+1)  & = & y_n(t) + \frac{d_n(t) - d_{n-1}(t)}{3} \label{ygapintermsofd}
\end{eqnarray} and replacing $t$ by $t+1$ in Eq. (\ref{dintermsofy}) and then plugging in Eq. (\ref{ygapintermsofd})
\aoj{
\begin{eqnarray*} d_0(t+1) - d_0(t) & = &  \frac{2(d_1(t)-d_0(t))}{3} \\
d_1(t+1) - d_1(t) & = &  \frac{d_{2}(t) - d_{1}(t)}{2} - \frac{d_1(t)-d_0(t)}{3} \\
d_i(t+1) - d_i(t) & = & \frac{d_{i+1}(t) - d_{i}(t)}{2}  - \frac{d_i(t) - d_{i-1}(t)}{2},~~~~i = 2, \ldots, n-1 \\
d_{n-1}(t+1) - d_{n-1}(t) & = & \frac{d_n(t) - d_{n-1}(t)}{3} - \frac{d_{n-1}(t) - d_{n-2}(t)}{2}\\
d_n(t+1) - d_n(t) & = & - \frac{2(d_n(t) - d_{n-1}(t))}{3} \\
\end{eqnarray*} or \[ d(t+1) - d(t) = \frac{1}{6} U_{n+1} d(t). \] {\bf q.e.d.}}
 
 \bigskip
 
 \noindent {\bf Lemma 5:} Let $k \geq 3$ and let $P_{k}=I+U_{k}/6$. Then the spectrum of $P_k$ is real. Labeling it from smallest to largest as
 $\lambda_{k}(P) \leq \cdots \leq \lambda_2(P) \leq \lambda_1(P)=1$, we have \[ \max( |\lambda_{k}(P)|, |\lambda_2(P)|) 
 \leq 1 - \frac{1}{\ao{3}k^2}. \]
 
 \noindent {\bf Proof:} For $k=3,4$ the lemma follows by an easy calculation. For general $k \geq 5$, consider an undirected line graph on $k$ nodes with self loops at each node. We assign weights to the edges as illustrated in Figure 1. Namely, the self loops on nodes $1,2,k-1,$ and $k$ have weight $1$; all other
self-loops have weight zero. The edge between nodes $1$ and $2$, and the edge between node $k-1$ and
$k$ each has weight $2$. Every other edge has weight $3$. We denote by $w_{ij}$  the weight of the edge
between $i$ and $j$ (defined to be zero if there is no edge between $i$ and $j$); since the graph is undirected, $w_{ij}=w_{ji}$ by definition. 
Moreover, we define $w_i$ to be  the sum of all the weights incident on 
node $i$, i.e. $w_i = \sum_{j=1}^n w_{ij}$. Thus, $(w_1,w_2,\ldots, w_{n}, w_{n+1}) = (3, 6, \ldots, 6, 3)$. With these definitions
in place, we observe that for all $i,j$, $P_{ij} = w_{ij}/w_i$. 

 \begin{figure}[h] \label{line}
\includegraphics[width=0.99\textwidth]{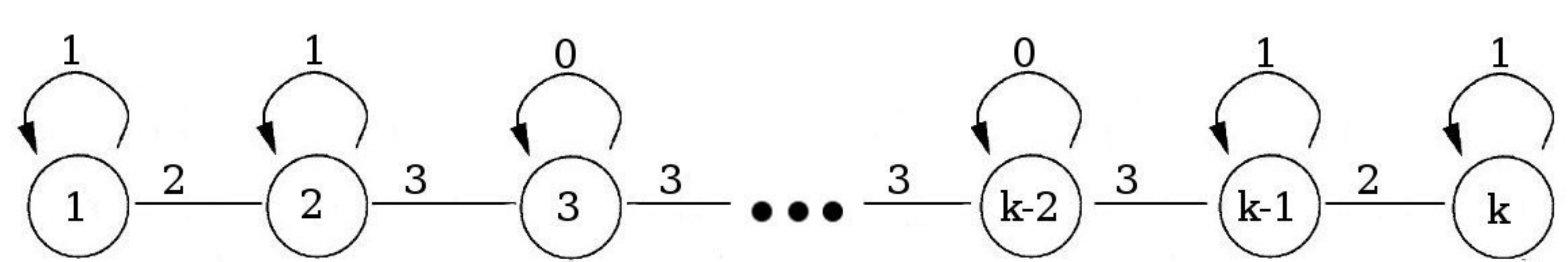}
\caption{The weighted graph capturing the transition matrix $P_k$ for $k \geq 5$: $P_{ij}=w_{ij}/w_{i}$ where $w_{ij}$ is the weight
of the edge $(i,j)$ and $w_i$ is the sum of all the weights incident on node $i$.}
\end{figure}

Define the inner product $\langle x, y \rangle = \sum_{i=1}^k w_i x_i y_i$. Then, \[ \langle x, Py \rangle = \sum_{i,j=1}^n w_{ij} x_i y_j = \sum_{i,j=1}^n w_{ji} x_i y_j = \langle Px, y \rangle, \] so that $P$ is self-adjoint in this inner product, and in particular its 
spectrum must be real.   Observing that the largest eigenvalue of $P$ is $1$ with the eigenvector of all ones,
the Courant-Fisher theorem gives \begin{small} \begin{eqnarray*} \lambda_2(P) & = & \max_{\langle x, \1 \rangle = 0, \langle x, x \rangle = 1} \langle x, Px \rangle  \\
& = & \max_{\sum_i w_i x_i = 0, \sum_i w_i x_i^2 = 1} \sum_{i,j} w_{ij} x_i x_j \\
& = & \max_{\sum_i w_i x_i = 0, \sum_i w_i x_i^2 = 1} \sum_{i,j} \frac{w_{ij}}{2} (x_i^2+ x_j^2 - (x_i - x_j)^2)  \\
& = & 1 - \min_{\sum_i w_i x_i = 0, \sum_i w_i x_i^2=1} \sum_{i<j} w_{ij} (x_i - x_j)^2 \\
& = & 1 - \min_{\sum_i w_i x_i = 0, \sum_i w_i x_i^2=1} \sum_{i=1}^{k-1} w_{i,i+1} (x_{i+1} - x_i)^2.
\end{eqnarray*} \end{small} We now lower bound the last term on the right hand side using a variation of the argument from $\cite{lo81}$. 
The minimum is achieved (since we are minimizing a continuous function over a compact set); use the notation $y$
for the minimizer. Consider the vector with $i$'th entry  $w_i y_i$. \aoc{Let $M$ be the index which maximizes $w_i y_i^2$; without loss of generality, we may assume that $y_M > 0$ (if not, we can replace $y$ by $-y$). Let $m$ be the index of the smallest entry; we may also assume without loss of generality that $m<M$.} Observe that the constraint $\sum_i w_i y_i = 0$ implies that $y_{\rm m}<0$ while the constraint $\sum_i w_i y_i^2 = 1$ implies
that the average value of $w_i y_i^2$ is $1/k$. This means $$w_M y_M^2 \geq 1/k$$ or $$y_M \geq 1/\sqrt{w_M k} \geq 1/\sqrt{6k}.$$ Thus
\[ \frac{1}{\sqrt{6k}} \leq y_M - y_m = \sum_{i=m}^{M-1} y_{i+1} - y_i.\] Applying Cauchy-Schwarz, 
\[ \frac{1}{6 k} \leq k \sum_{i=m}^{M-1} (y_{i+1}-y_i)^2, \] 
and therefore, 
\[ \sum_{i=1}^{k-1} w_{i,i+1} (y_{i+1}-y_i)^2  \geq \aoc{2 \sum_{i=m}^{M-1} (y_{i+1}-y_i)^2 \geq \frac{2}{6k^2} \geq }\frac{1}{\aoc{3}k^2}. \] Putting it all together, this
implies the desired bound on $\lambda_2$. 

A similar argument (\ao{which is also a variation of an argument from \cite{lo81}}) proves the bound for $\lambda_k$. By the Courant-Fisher theorem \begin{small}
\begin{eqnarray*} \lambda_{k}(P) & = & \min_{\langle x, x \rangle =1} \langle x, Px \rangle \\
& = & \min_{\sum_i w_i x_i^2=1} \sum_{i,j=1}^n w_{ij} x_i x_j \\
& = & \min_{\sum_i w_i x_i^2=1}  \sum_{i,j} \frac{w_{ij}}{2} ((x_i + x_j)^2 - x_i^2 - x_j^2 ) \\
& = & -1 + \min_{\sum_i w_i x_i^2 = 1} \sum_{i < j} w_{ij} (x_i + x_j)^2  + \sum_{i=1}^k 2 w_{i,i} x_i^2 \\
& = & -1 +  \min_{\sum_i w_i x_i^2 = 1} \sum_{i=1}^{k-1} w_{i,i+1} (x_{i+1}-x_i)^2 \\
&& ~~~~~~~~~~~~~~~~~~+ 2x_1^2 + \aoc{2} x_2^2 + \aoc{2} x_{k-1}^2 + 2x_{k}^2.
\end{eqnarray*} \end{small}where in the final step we can replace $(x_{i+1} \aoc{+} x_i)^2$ with $(x_{i+1} \aoc{-} x_i)^2$ since 
that flipping the sign of every other $x_i$ does not change the constraint $\sum_i w_i x_i^2$.  Next, as before let $y$ be the vector achieving 
the optimum in the above expression, and let $w_M y_M$ be the largest among all the numbers $w_i y_i$. As before, $y_M \geq 1/\sqrt{6k}$. If
$M=1$ or $M=2$, the above equation immediately proves the lemma.  If $M>2$ \[ \frac{1}{\sqrt{6k}} \leq y_{\aoc{M}} = y_2 + (y_3-y_2) + \cdots + (y_M - y_{M-1}), \] so that by Cauchy-Schwarz, 
\[ \frac{1}{6k} \leq k(y_2^2 + \sum_{k=2}^{n-1} (y_{k+1} - y_k)^2), \] \aoc{so that \[ \min_{\sum_i w_i x_i^2 = 1} \sum_{i=1}^{k-1} w_{i,i+1} (x_{i+1}-x_i)^2 + 2x_1^2 + \aoc{2} x_2^2 + \aoc{2} x_{k-1}^2 + 2x_{k}^2 \geq 2 \min_{\sum_i w_i x_i^2 = 1} x_2^2 + \sum_{i=2}^{k-1} (x_{i+1}-x_i)^2 \geq \frac{1}{3k^2}, \] }which
implies the corresponding bound on $\lambda_{k}$. {\bf q.e.d.} 

\smallskip

\noindent {\bf Proof of Theorem 1:}  Decompose $d(0)$ in terms of the eigenvectors of $P$ as
\[ d(0) = \sum_{i=1}^{n+1} c_i v_i, \] where $v_i$ \aoc{corresponds to eigenvalue} $\lambda_i$. We know that 
$v_1 = \1$. Thus  %and $c_1 = \langle d(0), \1 \rangle/\langle 1, 1 \rangle = (1/6n) \sum_{i=1}^{n+1} w_i d_i(0)$. Moreover
\[ d(t) = c_1 \1  + \sum_{i=2}^{n+1} \lambda_i^t c_i v_i.\] Using the fact that the eigenvectors $v_i$ are orthogonal 
in the inner product $\langle \cdot , \cdot \rangle$,
\begin{small} \[ \langle d(t) - c_1 \1, d(t) - c_1 \1 \rangle \leq (1-\frac{1}{3(n+1)^2})^{\aoc{2t}} 
\langle d(0) - c_1 \1, d(0) - c_1 \1 \rangle \] \end{small} or 
\[ \sum_i (d_i(t) - c_1)^2 \leq \frac{1}{3} (1-\frac{1}{3(n+1)^2})^{\aoc{2t}} \langle d(0) - c_1 \1, d(0) - c_1 \1 \rangle.\] Observe that 
%\[  \langle d(0) - c_1 \1, d(0) - c_1 \1 \rangle  \leq  \langle d(0) , d(0) \rangle  \ao{\leq} 6 F(1) \]
%\aoj{\begin{eqnarray*} \langle d(0) - c_1 \1, d(0) - c_1 \1 \rangle & \leq & \langle d(0) , d(0) \rangle \\ 
%& \leq & 6 \sum_{i=1}^{n+1} d_i^2(0) \\
%& \leq & 6 \sum_{i=1}^{n+1} d_i(0) \\
%& \leq & 6 F(1), \end{eqnarray*}  }
\aoj{\[  \langle d(0) - c_1 \1, d(0) - c_1 \1 \rangle  \leq  \langle d(0) , d(0) \rangle 
 \leq  6 \sum_{i=1}^{n+1} d_i^2(0) 
 \leq  6 (\sum_{i=1}^{n+1} d_i(0))^2
 \leq  6 F(1)^2 \]  }
so that 
\[ \sum_i (d_i(t)-c_1)^2 \leq 2 (1-\frac{1}{3(n+1)^2})^{\aoc{2t}} F(1)^2.\] It follows that each $d_i(t)$ has a limit \aoc{(which equals $c_1$)}, and this immediately implies that
every $y_i(t)$ has a limit, and therefore every $x_i(t)$ has a limit. Most importantly,
because every $d_i(t)$ approaches the same limit, we can conclude that the optimality condition of Lemma 2 is satisfied in the limit, and the set of 
limiting positions achieves $\Phi^*$. 

Moreover, after $\aoc{\frac{3}{2}}(n+1)^2 \log (2F(1)^2n^2/\epsilon^2)$ iterations, each
$d_i(t)$ is within $\epsilon/n$ of its final value, which implies every $y_i(t)$ is within $\epsilon$ of its final value. 

Now the static control law is invariant under scaling of $\rho$, so we can assume without loss of generality that we are dealing with $\widehat{\rho} = \rho/\rho_{\rm min}$, so that the minimum value 
of $\widehat{\rho}$ is $1$. Then, we conclude that after $3(n+1)^2 \log ( \sqrt{2}  \int_0^1 \hat{\rho}(z) dz \;  n/\epsilon)$ iterations, every $x_i(t)$ is within $\epsilon$ of its final value. To conclude the proof,
observe that $\int_0^1 \widehat{\rho}(z) dz \leq \rho_{\rm max}/\rho_{ \min }$, and the upper bound in the theorem statement follows.  {\bf q.e.d}.

\smallskip

\noindent {\bf Remark:} Observe that the coverage control law we have presented in this section is naturally robust to addition and deletion of 
agents as well as changes in $\rho$. Indeed, as long as the density and the number of agents stop changing at some point, the algorithm is
guaranteed to converge to the optimal configuration. An open problem is to prove performance bounds for this algorithm in the scenario when the number of agents and the density are continually changing.  

%%%%%%%%%%%%%%%%%%%%%%%%%%%%%%%%%%%%%%%%%%%%%%%%%%%%%%%%%%%%%%%%%%%%%%%%%%%%%%%%
\section{A DYNAMIC COVERAGE CONTROL LAW} In this section, we propose another control law for the nonuniform coverage problem on the line. 
\aoc{Our main result is that this control law converges an order of magnitude better than the static control law of the previous section; in particular, it drives the agents \aoc{close} to the optimal positions in a number of \aoc{communication/movement rounds} that is
\aoc{essentially} linear in $n$ for each agent. This scaling is optimal in the sense that any control law requires at least linearly many discrete-time rounds of communication to drive agents close to their optimal positions.}

We draw heavily on the paper \cite{dhn00}, which described a fast ``nonreversible Markov sampler'' for sampling  a uniform random 
number from $\{1,\ldots,n\}$. \aoc{The main technical component of our result is a modification of the results of \cite{dhn00}; specifically, we use the ideas of \cite{dhn00} to modify the nonreversible Markov sampler to handle the nonuniform sampling problem of coverage control.}

We make stronger assumptions on the capabilities and knowledge of each agent. We now assume that
agents can store numbers in memory, transmit numbers to their neighbors, and can detect when their neighbors move to a new location.  \ao{However, every node only stores and transmits/senses two numbers in each of the control
law iterations, so  the additional effort expended is not excessively onerous.} Moreover, we assume that each agent has an estimate $U$ of the total number of agents, and that this estimate is accurate
within a constant factor $c$ \aor{which does not depend on $n$}:  \begin{equation} \label{uassumption} \frac{n}{c} \leq U \leq c n. \end{equation}\aoj{Finally, we assume
for simplicity that the number of agents is at least $3$, and consequently $U \geq 3$ as well.} 

We refer to the control law of this section as
the ``dynamic coverage control law,'' since in contrast with the control law of the previous section, the feedback law has 
dynamics of its own. This control law follows two steps: an initial measurement step and a subsequent communication/measurement/movement 
step.

  \begin{figure*}[t]
\centering \label{lifted}
\includegraphics[width=1.05\textwidth]{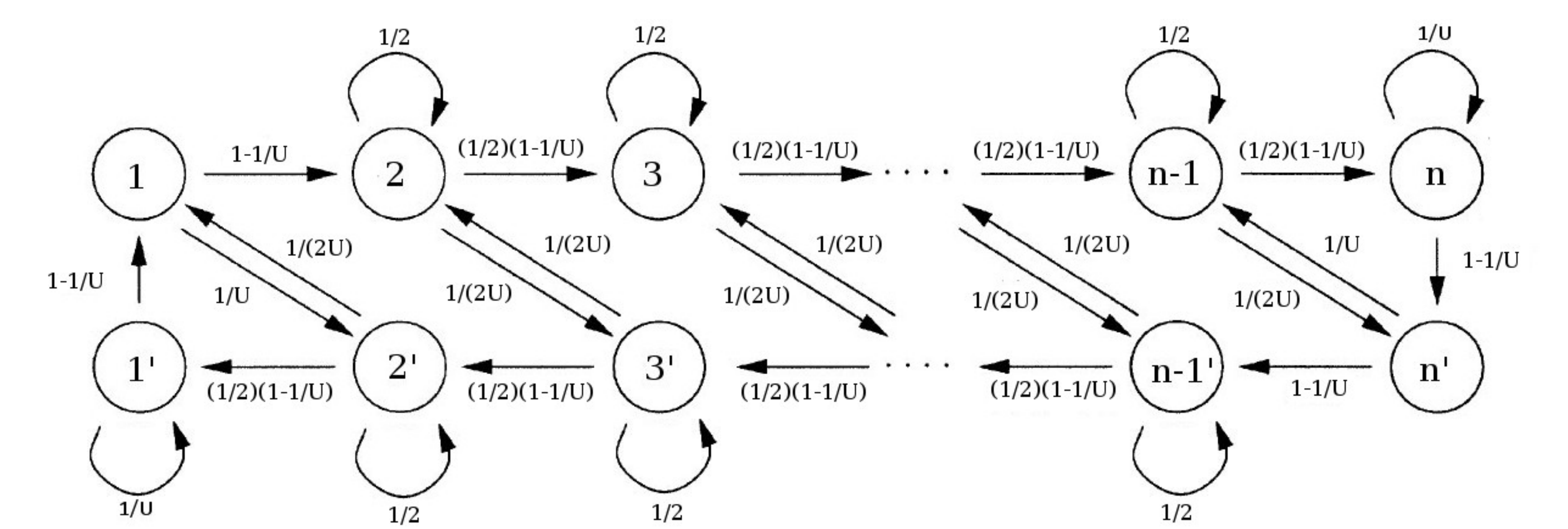}
\caption{The Markov chain representing the update equations of the dynamic coverage control law.}
\end{figure*}

\smallskip

\noindent {\bf Dynamic control law:} Nodes keep track of the variables $z_i(t), z_{i'}(t)$, initialized in the first step as 
\begin{eqnarray*} z_1(0) & = & \frac{1}{2}  \int_0^{m_{\rho}^1(x_1(0), x_2(0))} \rho(z) dz \\
z_i(0) & = & \frac{1}{2} \int_{m_{\rho}^1 (x_{i-1}(0), x_{i}(0))}^{m_{\rho}^1(x_i(0), x_{i+1}(0))} \rho(z)dz ~~~i = 2,\ldots,n-1 \\
z_n(0) & = & \frac{1}{2} \int_{m_{\rho}^1(x_{n-1}(0), x_n(0))}^1 \rho(z) dz 
\end{eqnarray*} and $z_{i'}(0)=z_i(0)$ for each $i$. At each step, nodes transmit their variables $z_i(t), z_{i'}(t)$ to their
neighbors, and then set their values $z_{i}(t+1), z_{i'}(t+1)$ to be linear combinations of their previous values and
the values they have just received. The linear combination taken by each agent is derived from a rule
based on Figure 2. Note that this figure contains nodes labeled $1,\ldots,n$ 
and $1',\ldots,n'$. Node $i$ sets $z_i(t+1)$ to be a linear combination of those values $z_{k}(t)$ which have edges going from $k$
to $i$; the coefficient it puts in front of $z_k(t)$ is the label on the edge. The value of $z_{i'}(t+1)$ is determined 
in the same way.  For example, agent $1$ updates as
\begin{eqnarray*} z_1(t+1) & = & (1 - \frac{1}{U}) z_{1'}(t) + \frac{1}{2U} z_{2'}(t) \\
z_1'(t+1) & = & \frac{1}{2} (1 - \frac{1}{U}) z_{2'}(t) + \frac{1}{U} z_{1'}(t). \\
% & & \\
%z_2(t+1) & = & (1 - \frac{1}{U}) z_1(t) + \frac{1}{2U} z_3'(t) + \frac{1}{2} z_2(t) \\
%z_2'(t+1) & = & \frac{1}{2} (1 - \frac{1}{U}) z_3'(t+1) + \frac{1}{U} z_1(t) + \frac{1}{2} z_2'(t) \\
%& & \\
%z_{n-1}(t+1) & = & \frac{1}{2} (1 - \frac{1}{U}) z_{n-2}(t) + \frac{1}{U} z_{n}'(t) + \frac{1}{2} z_{n-1}(t) \\
%z_{n-1}'(t+1) & = & ( 1 - \frac{1}{U}) z_n'(t) + \frac{1}{2U} z_{n-2}(t) + \frac{1}{2} z_{n-1}'(t) \\
%& & \\
%_n(t+1) & = & \frac{1}{2} ( 1 - \frac{1}{U}) z_{n-1}(t) + \frac{1}{U} z_n(t) \\
%z_n'(t+1) & = & ( 1 - \frac{1}{U}) z_n(t) + \frac{1}{2U} z_{n-1}(t)
\end{eqnarray*} %while the remaining agents $3,\ldots, n-2$ update as
%\begin{eqnarray*}
%z_i(t+1) & = & \frac{1}{2} (1 - \frac{1}{U}) z_{i-1}(t) + \frac{1}{2U} z_{i+1}' (t) + \frac{1}{2} z_i(t) \\
%z_i'(t+1) & = &  \frac{1}{2} (1 - \frac{1}{U}) z_{i+1}'(t) + \frac{1}{2U} z_{i-1}(t) + \frac{1}{2} z_i'(t) \\
%\end{eqnarray*}
After updating their variables $z_i, z_{i'}$, \aor{some of the agents execute motions according to the following rule. Initially, at round $1$ the leftmost agent  initiates the motion by moving to the 
point $c$ satisfying \[ \int_0^c \rho(z) ~dz = z_1(1) + z_{1'}(1). \] If this number $c$ is higher than $x_2(1)$, then agent $2$ also moves to point $c$, 
as well as all other agents whose positions at time $1$ are less than $c$. This may be implemented in the straightforward way: when agent $1$ reaches
the location of agent $2$, the two begin moving together. }

\aor{Generally, after agent $i$ initiates the motion at round $t$,  agent $i+1$ initiates the motion at time $t+1$. It moves to the position $c$ which 
satisfies \[ \int_{x_{i}(t)}^{c} \rho(z) dz = z_i(t+1) + z_{i'}(t+1),\] and every agent whose position is strictly less than $c$ moves to position $c$ alongside it. 
Every $U$ steps, the leftmost agent initiates the motion once again and the process repeats itself. }
%\begin{center} \begin{figure}[h] 
%        \epsfig{file=lifted.eps, scale=0.2} 
%  \end{figure} \end{center}

\bigskip

\aor{\noindent {\bf Remark:} We do not consider the issue of collision avoidance in this paper, and thus the description of the algorithm above allows two 
agents to share the same location. It is possible to modify the dynamic control law so that no two agents are ever in the same place by introducing a tolerance $\epsilon$ so that node $i$ begins moving whenever node $i-1$ is within $\epsilon$ distance of it. For reasons of simplicity we will not be considering these issues
here. }

\bigskip

While the \aor{update rule of the dynamic control law} is somewhat involved, it has a simple interpretation. Consider the Markov chain
of Figure 2 where we interpret the label
on the edge $(i,j)$ as the probability of transitioning from $i$ to $j$; note that the labels on the outgoing edges sum to $1$ for every node. If $z_i(t)$ is the probability of being at node $i$ at time $t$, and
$z_i'(t)$ is the probability of being at node $i'$ at time $t$, then the variables $z_i(t), z_{i'}(t)$ satisfy
the above recursion. 

Indeed, this recursion is an adaptation of the ``nonreversible Markov chain sampler'' from \cite{dhn00}. It was observed in 
that paper that while an ordinary ``diffusive'' Markov chain on the line graph which, say, moves to the right and left each with probability $1/2$ takes on
the order of $n^2$ steps to come close to the uniform distribution, \aoc{a  ``guided'' Markov chain which moves clockwise round the circle with probability $1-1/n$ and transitions to a symmetrically placed node with probability $1/n$ gets close to the uniform distribution an 
order of magnitude faster. } The dynamic coverage control law of this section is an attempt to harness this insight for the
coverage problem. 

Our goal in this section is to prove that this iteration solves the coverage control problem and to demonstrate that it is an order of 
magnitude \aor{better as far as number of update rounds} than the static control law of the previous 
section. The main result of this section is the following theorem. 

\bigskip

\noindent {\bf Theorem 6:} Each $x_i(t)$ has a limit, and the limiting set of positions has coverage $\Phi^*$. Moreover, after
every agent has executed $O(n \log ( \frac{ \ao{\rho_{\rm max}} n}{\ao{\rho_{\rm min}} \epsilon} ))$ rounds of updates, each agent is within $\epsilon$ of its final 
limit. 

\bigskip

\aor{\noindent {\bf Remark:} We remark that while this is an order of magnitude improvement in number of update rounds, understood here as a step of sensing, communication, computation, and movement executed in parallel by all the nodes, it is not necessarily an improvement in speed or energy. Indeed, an update round in the dynamic control law includes communication
between neighboring agents, while an update round in the static control law does not. Moreover, the time it takes to execute the motion in each round in the two algorithms is not immediately comparable; note that, assuming that each node can travel at unit speed, executing the motions takes $\max_{i=1}^n x_i(t+1) - x_i(t)$ time in the first algorithm and $x_j(t+1) - x_j(t)$ in the second, where $j$ is the index initiating movement at time $t$. However, since the two algorithms send agents to different positions, the movement times are not immediately comparable.  Whether the dynamic algorithm is in fact faster in time or takes less energy will depend on the precise speed and energy requirements of 
sensing, communication, computation, and movement. Our main results, in Theorems 1 and 6, only demonstrate the weaker statement that the number of rounds 
in the dynamic algorithm is an order of magnitude better.}

\aor{Nevertheless, we note that under some additional assumptions, an improvement in the number of  rounds does translate into an improvement in speed and energy. For example, if the time taken and energy expended are dominated by sensing and communication, and
if the acts of sensing the field $\rho$ at each step and exchanging messages with neighbors take
similar expenditures of time and energy, then an order of magnitude improvement in the number of rounds
immediately translates into an order of magnitude improvement in speed
and energy.}

%\bigskip
%
%\aor{\noindent {\bf Remark:} We note, however, that under a number of plausible assumptions on relative costs of sensing, communication, computation, and movement, the dynamic control law is indeed an order of magnitude better in both time taken and energy expended. We briefly outline one such set of assumptions next. Suppose that communication and sensing take roughly similar time and energy expenditure, and computation is negligible. A consequence of Theorem 6 is 
%that the amount of movement executed by each agent approaches zero as time goes to infinity. Indeed, for every $\epsilon$, there is a time when the agent is
%always within $\epsilon$ of its final position, so its movements at each subsequent round cannot be bigger than $2 \epsilon$. It follows that, in the long run, the cost of movement 
%end up being negligible. Consequently, as long communication among neighbors takes roughly the same time/energy as sensing the field $\rho$, and expenditures on computation are negligible compared to both, an improvement in the number of rounds translates into an improvement in both time taken and energy expended as $t \rightarrow \infty$. To put it another way, if we consider the time taken or energy expended until every node is within $\epsilon$ of its final position under these assumptions, then as $\epsilon \rightarrow 0$, the time taken and energy expended by the dynamic control law is an order of magnitude smaller.} 

\bigskip

The \aoc{general structure} of our proof is similar the proof of the related results from \cite{dhn00}. Specifically, we prove 
the theorem by first establishing a certain inequality on the state of the Markov chain of Figure 2, namely, Eq. (\ref{spreading}); the proof of \cite{dhn00} proceeds in the same way. \aor{The meaning of Eq. (\ref{spreading}) is that after $O(n)$ steps, the state of the chain is not too far from uniform in a certain sense. Due to the  somewhat asymmetric nature of our chain, our proof of this inequality is considerably more involved than in \cite{dhn00} and occupies the bulk of this section. Thus we begin with an extended series of lemmas whose goal is to prove Eq. (\ref{spreading}); these lemmas essentially describe the action of the probability transition matrix in a slightly different way which is more convenient for the proof of Eq. (\ref{spreading}). Once Eq. (\ref{spreading}) is proved, Theorem 6 will follow by a
standard argument in probability.}

\smallskip

\aoc{We now proceed with our proof.} Let  $z(t)$ denote the row vector 
$\left( \begin{array}{cccccc} z_1(t) & \ldots & z_n(t) & z_{1'}(t) & \ldots z_{n'}(t) \end{array} \right)$. 
Let $K$ denote the matrix which maps $z(t)$ to $z(t+1)$ through right-multiplication: 
\[ z(t+1) = z(t)  K. \] Observe that $K$ is a nonnegative, irreducible stochastic matrix. Standard results in
Markov chain theory imply that the above iteration converges to a scaled multiple of the stationary probability
vector of the chain. Thus we can immediately conclude that each $z_i(t)$ has a limit, and consequently, 
the positions of the agents under the dynamic coverage control law have limits as well. 

%Moreover, observe that the
%vector of all ones is a right eigenvector of $K$, which implies that $\sum_{i=1}^n z_i(t) + z_{i'}(t)$ does not change after the %execution of each update. Since $\sum_{i=1}^n z_i(0) + z_{i'}(0)=F(1)$, we can immediately conclude that no agent ever moves outside of %$[0,1]$. 

\smallskip

\noindent {\bf Lemma 7:} The stationary probability of the Markov chain in Figure 2 is
\[ \pi_1 = \pi_{1'} = \pi_n = \pi_{n'} = \frac{1}{4 \aoc{n}} \] and for all other nodes, \[ \pi_i = \pi_{i'} = \frac{1}{2\aoc{n}}. \]

\bigskip 

This lemma can be proved by verifying that the stacked vector $\pi$ satisfies $\pi^T K = \pi^T$. 

This lemma implies that \begin{equation} \label{1limit} \lim_{t} z_1(t) = \lim_{t} z_{1'}(t) = \lim_{t} z_n(t) = \lim_{t} z_{n'}(t) = \frac{F(1)}{4\aoc{n}}, \end{equation} 
and \begin{equation} \label{ilimit} \lim_{t} z_i(t) = \lim_{t} z_{i'}(t) = \frac{F(1)}{2 \aoc{n}}, \end{equation} which implies by Lemma 2 that the limiting set of positions 
do have optimal coverage $\Phi^*$. It remains to bound the time until the agents approach these positions. \aoj{We begin with several definitions
and lemmas which we make use of later.}

\bigskip

\aoj{\noindent {\bf Definition:} Let ${\cal A}(n_1, \ldots, n_k)$ be the set of sequences of symbols $A_1, A_2, \ldots, A_k$ in which the symbol $A_k$ appears $n_k$ times.}

\bigskip

\aoj{\noindent {\bf Lemma 8:} Consider the following process of randomly generating elements from ${\cal A}(n_1, \ldots, n_k)$. Begin with the sequence of $A_1$ repeated $n_1$ times; then place an $A_2$ at a uniformly random slot, where a slot is either between the letters already in place, or before the first or after the last letter; repeat the process $n_2$ times. Finally, perform the analogous operations for $A_3, \ldots, A_k$, that is, place an $A_3$ at a uniformly random slot $n_3$ times, place an $A_4$ at a uniformly random slot $n_4$ times, etc. Then every sequence in ${\cal A}(n_1, \ldots, n_k)$ has equal probability of being generated by this process.} 

\bigskip

\aoj{\noindent {\bf Example:} The proof of this lemma will be made more intuitive by first considering an example. Consider the sequence $A_1 A_2 A_2 A_1 A_2 \in {\cal A}(2,3)$. It may be generated with this process in exactly six sequences of insertions, shown below. The number before each sequence is the slot (enumerated from left to right) into which the next letter will be inserted.   \begin{small} $$\begin{array}{cccccc} 
2.A_1 A_1                &  2.A_1 A_1                &  2.A_1  A_1               &  2.A_1 A_1                &  3.A_1 A_1                &  3.A_1 A_1 \\
2. A_1 A_2 A_1         &  3. A_1 A_2 A_1         &  4. A_1 A_2 A_1         &  4. A_1 A_2 A_1         &  2. A_1 A_1 A_2         &  2. A_1 A_1 A_2 \\
5. A_1 A_2 A_2 A_1     &  5. A_1 A_2 A_2 A_1     &  2. A_1 A_2 A_1 A_2     &  3. A_1 A_2 A_1 A_2     &  2. A_1 A_2 A_1 A_2     &  3. A_1 A_2 A_1 A_2\\
A_1 A_2 A_2 A_1 A_2 &   A_1 A_2 A_2 A_1 A_2 &  A_1 A_2 A_2 A_1 A_2 &   A_1 A_2 A_2 A_1 A_2 &   A_1 A_2 A_2 A_1 A_2 &  A_1 A_2 A_2 A_1 A_2 \\
\end{array}$$ \end{small}Note that each of the above sequences has probability $(1/3) (1/4) (1/5)$ of being generated by the process of Lemma 8.}

\bigskip

\aoj{\noindent {\bf Proof:} Observe that the insertion of a symbol into a slot replaces it with two slots, so that every choice of slots has the same probability. It therefore suffices to show that any sequence can be generated in the same number of ways using the process described in this lemma.}

\aoj{Consider the reverse process, i.e., the sequence of deletions beginning with the final sequence - first of $A_k$'s, then of $A_{k-1}$'s, and so on - until we are left with the sequence of $n_1$ $A_1$'s. For example, the possible deletions of the sequence $A_1 A_2 A_2 A_1 A_2$ are
\begin{small} $$\begin{array}{cccccc} 
 5. A_1 A_2 A_2 A_1 A_2 &   5. A_1 A_2 A_2 A_1 A_2 &   2.A_1 A_2 A_2 A_1 A_2 &  3. A_1 A_2 A_2 A_1 A_2 &  2. A_1 A_2 A_2 A_1 A_2 &   3. A_1 A_2 A_2 A_1 A_2 \\
2. A_1 A_2 A_2 A_1      &  3. A_1 A_2 A_2 A_1     &   4.A_1  A_2 A_1 A_2     &   4. A_1  A_2   A_1 A_2     &  2. A_1  A_2 A_1 A_2     &   2. A_1 A_2  A_1 \\
2. A_1  A_2 A_1          &   2. A_1  A_2  A_1          &   2.A_1   A_2 A_1          &  2. A_1 A_2   A_1         &   3. A_1  A_1 A_2         &   3. A_1  A_1 A_2 \\
A_1  A_1                 &  A_1   A_1                 &  A_1   A_1                &  A_1   A_1                &  A_1 A_1                 &  A_1 A_1  \\
\end{array}$$ \end{small}where the number next to the sequence is the index of the letter being deleted. As the above example illustrates, every sequence of insertions naturally corresponds to a sequence of deletions: if one can obtain sequence $s_2$ by insertion of a letter into slot $j$ of sequence $s_1$, then one can obtain $s_1$ by deleting the $j$'th letter from $s_2$. Thus, for any sequence, there are as many ways to generate it through the insertion process of this lemma as there to generate the sequence of $A_1$ repeated $n_1$ times by this deletion process, and there are plainly $n_{k}! n_{k-1}! \cdots n_2!$ ways to do the latter. Since the latter number does not depend on the particular sequence, this completes the proof. {\bf q.e.d.}}

\bigskip

\aoj{\noindent {\bf Definition:} Consider the Markov chain obtained by adding a self-loop of probability $1/2$ to nodes $1,1',n,n'$ in Figure 2 and scaling the probabilities on the outgoing links from those nodes by a factor of $1/2$. We  refer to it as ``the modified Markov chain'' to distinguish it from the 
Markov chain of Figure 2, to which we henceforth refer as ``the original Markov chain.'' The modified Markov chain is appropriately viewed as a certain 
random walk on a multigraph, as nodes $1'$ and $n$ have two self-loops.} 

\bigskip

\aoj{\noindent {\bf Definition:} Observe that every node in the modified Markov chain has three outgoing edges: a self-loop with probability $1/2$, a link with probability
$1/2-1/(2U)$, and a link with probability $1/(2U)$. We refer to the first as ``major self loops,'' to the second as ``continuing links,'' and to the third as
``switching links.''} \aoc{Note that the switching links originating at nodes $1'$ and $n$ are self-loops.} 

\bigskip

\aoj{\noindent {\bf Lemma 9.} Let $T(a,j)$ be the event that starting at node $a$ at time $t=1$, the modified Markov chain satisfies the following conditions \aoc{during times\footnote{\aoc{Henceforth, we will shorten this statement as ``by time $4n$.''}} $1,\ldots, 4n-1$}: 
\begin{enumerate} 
\item The chain has not taken major self-loops at nodes $1,1',n,n'$ . 
\item A single switching link was taken, at node $j$. 
\end{enumerate} Define $L_k$ be the event that $k$ major self-loops were taken by time $t=4n$.
Then there is some constant $c'$ which does not depend on $n$ and $U$ such that \[ P(T(a,j) ~|~ L_k) \geq \frac{1}{c'n} ~~~~ \mbox{ for all } a,j \mbox{ and } k \leq 2n-2.\]}

\noindent {\bf Proof:} \aoc{We first outline the general structure of the proof, which proceeds by repeated conditioning. Instead of conditioning on $L_k$, we will instead condition on the specific set of $k$ times when the modified chain took a self-loop (we call this $L(4n)$); specifically, we argue for any $L(4n)$ there is a constant probability that exactly one switching link was taken by time $4n$. We then argue that conditional on $L(4n)$ {\em and} a single switching link being taken by time $4n$, the probability that that link was taken at node $j$ and that no self-loops were taken at nodes $1,1',n,n'$ is proportional to $1/n$ (this last step crucially relies on Lemma 8). This proves the current lemma. }

\aoc{We now proceed with the proof. As just mentioned, we define  $L(4n)$ to be} the set of times in $\{1,\ldots,4n-1\}$ when the modified Markov chain takes a major self-loop. Because the distribution of the number of major self-loops taken by the modified chain is exactly that of the times a fair coin lands on heads in $4n-1$ tosses, we have that conditional on $L_k$, $L(4n)$ is uniformly distributed over all $k$-element subsets of $\{1,\ldots,4n-1\}$. Similarly, conditional on any $k$-element $L(4n)$, the distribution of the number of switching links is exactly that of the number of times a coin whose heads probability is $1/U$ lands on heads \aoc{exactly once}  in $4n-1-k$ tosses. So, if $S$ is the event that a single switching link was taken by time $t=4n$, then for any $k$-element $L(4n)$,\[ P(S|L(4n)) = (4n-1-k) \frac{1}{U} (1-\frac{1}{U})^{4n-1-k} \geq \frac{2n}{U} ((1 -\frac{1}{U})^U)^{4(n/U)} \geq 2c (\frac{2}{3})^{12c}, \] where we used the bound $k \leq 2n-2$ and that $(1-1/x)^x$ is increasing for $x \geq 1$ while $U \geq 3$. It follows that conditional on $L_k$, $S$ occurs with constant probability independent of $n,U$.

\aoj{Moreover, conditional on both $L_k$ and $S$, the sequence of continuing links, switching links, and major self-loops has uniform distribution over the set of sequences with $k$ major self-loops, a single switching link, and $4n-2-k$ continuing links. By Lemma 8 we  may generate this sequence by writing down 
$4n-2-k$ C's (corresponding to continuing links), inserting an S (corresponding to a switching link) in a uniformly random slot, and then inserting $k$ L's 
(corresponding to major self-loops) in uniformly random slots. To conclude the proof of the lemma, we next argue that for any node $j$, with probability $c''/n$, (where $c''$ is some constant independent of $n,U$) this process produces a sequence which takes a switching link at node $j$ and avoids self-loops at nodes $1,1',n,n'$. }

\aoj{Indeed, there are $4n-2-k \geq 2n$ continuing links (where we used that $k \leq 2n-2$), so every node is visited at least once, and the probability of inserting a switching link in the corresponding slot is at least $1/(4n)$. Moreover, because there are at most $4n$ continuing links and at most one switching link, the number of times the chain can enter the sets $\{n',1\}$ and $\{n,1'\}$ is at most $5$. Thus we avoid creating  a major self-loop at $1,1',n,n'$ if and only if we avoid at most ten slots every time we insert an $L$ into a random slot. When $n \geq 6$, the probability of this is at least
$$ (1 - \frac{10}{4n-1-k})^{4n} \geq (1-\frac{10}{2n})^{4n} \geq (\frac{1}{6})^{24.}$$ To summarize, for $n \geq 6$, the probability of placing a link at $j$ and avoiding major self-loops at nodes $1,1',n,n'$ conditional on $L_k,S$ is at least $(1/(4n)) (1/6)^{24}$. This proves the lemma for $n \geq 6$, and it holds automatically for $n=3,4,5$ (because for these $n$, the probability in question is positive, and so can be assumed to be lower bounded by a constant which can be taken to be independent of $n,U$).  {\bf q.e.d.}}

\bigskip

\aoj{\noindent {\bf Lemma 10:} Let $T(a,b)$ be the event that: 
\begin{enumerate}
\item[i.] The  modified chain starts at node $a$ at time $t=1$ and is at node $b$ at time $t=4n$. 
\item[ii.]  The modified chain has not taken the self loops at nodes $1,1',n,n'$ before time $t=4n$. 
\end{enumerate} If $k$ is at most $2n-2$ and has the same parity as $b-a$, then
there is some $j$ (depending on $k$) such that \[ T(a,j) \cap L_k \subset T(a,b) \cap L_k.\]   }

\aoj{\noindent {\bf Proof:} In words, this lemma says that conditional on having taken a certain number of self-loops which
is at most  $2n-2$ and has the right parity, $T(a,b)$ contains some $T(a,j)$. To prove this, label the nodes of the modified Markov chain $1, \ldots, 2n$ counterclockwise, and let $p(t)$ be the label of the state of the modified Markov chain at time $t$ taken
${\rm mod}~ 2n$. \aoc{This simply means that we will refer to state $2n$ as state zero, but the advantage of such a labeling is that  some facts become easier to state in terms of addition ${\rm mod}~ 2n$.}  In particular,  a major self loop does not change $p(t)$, a continuing link adds $1$ to $p(t)$, and a switching link taken at node $n + z$ adds $-2z$ to $p(t)$.}

\aoj{With that in mind, we argue that a chain beginning at node $a$ and by time $4n$ taking $k$ major self-loops, one switching link at node $n+z$, where $z$ is specially chosen, and continuing links for the remainder of the time will be located at $b$ at time $t=4n$. This will prove the lemma. }

\aoj{Indeed, putting off the choice of $z$ for the moment, we have that
\[ p(4n) = a+(4n-1-1-k) - 2z. \] Thus $p(4n)=b$ if and only if 
$$ 4n-2-k-2z = b-a.$$ Because $k$ has the same parity as $b-a$, this equation has a solution and so a satisfactory $z$ can always be found. {\bf q.e.d.}}

\bigskip

\aoj{\noindent {\bf Lemma 11:} \[ P(T(a,b)) \geq \frac{1}{c'' n} ~~~~\mbox{ for any } a,b,\] where
$c''>0$ is a constant which does not depend on $n$ and $U$. }

\bigskip

\aoj{\noindent {\bf Proof:} If $b-a$ is even,  \begin{eqnarray*} P(T(a,b)) & \geq & \sum_{k ~{\rm even}, ~~ k \leq 2n-2} P(k \mbox{ self loops taken}) P(T(a,b)~|~k \mbox{ self loops taken} ) \\
& \geq & \sum_{k ~{\rm even},~~ k \leq 2n-2} P(k \mbox{ self loops taken}) P(T(a,j)~|~k \mbox{ self loops taken} ) \\ 
& \geq & \frac{1}{c'n} \sum_{k ~{\rm even},~~ k \leq 2n-2} P(k \mbox{ self loops taken}) \\
& \geq &  \frac{1}{c''n}.
\end{eqnarray*} where the transition between the first and second line used Lemma 10 and the transition between the second and 
third used Lemma 9. The case when $b-a$ is odd is similar. {\bf q.e.d.}}

\bigskip

\aoj{With the last lemma in place, we are finally able to turn to the analysis of the original Markov chain from Figure 2.} We find it convenient to use the largest $l^1$ distance between the rows of $K^t$ and their
ultimate limit as a measure of convergence. In particular, let $(K^t)_i$ refer to the $i$'th row of $K^t$, and let 
\[ v(t) = \max_{i} ||(K^t)_i - \pi||_1. \] We use $t(n,\epsilon)$ to denote the time until $v(t)$ permanently
sinks below $\epsilon$.

\noindent {\bf Lemma 12} $t(n,\epsilon)=O(n \log 1/\epsilon)$. 

\smallskip

\aoj{ \noindent {\bf Proof:} By Lemma 11, the modified Markov chain starting from  node $a$ at time $t=1$ has a $1/(c''n)$ chance of 
being at any node $b$ at time $t=4n$ and avoiding the major self-loops at nodes $1,1',n,n'$ while getting there. But all choices of links satisfying these two conditions have a higher probability 
in the original Markov chain which has no major self-loops at $1,1',n,n'$. Consequently, if $s(t)$ is the state of the original Markov chain at time $t$, then \begin{equation}  P(s(4n)=b ~|~ s(1)=a) \geq \frac{1}{c'' n} ~~~~ \mbox{ for all } a,b,n, \label{spreading} \end{equation}  where $c''$ is the same constant as in Lemma 11. It is not hard to see (and was shown in \cite{dhn00}) that Eq. (\ref{spreading}) implies the conclusion of this lemma. {\bf q.e.d}}

\noindent {\bf Lemma 13} After $O(n \log  (\ao{nF(1)} /\epsilon) )$ steps, we will always have 
%\[ |z_1(t) - \frac{S}{4n}| \leq \frac{\epsilon}{2n}, \mbox{    and    } |z_{1'}(t) - \frac{S}{4n}| \leq \frac{\epsilon}{2n} \] 
%\[ |z_n(t) - \frac{S}{4n}| \leq \frac{\epsilon}{2n}, \mbox{    and    } |z_{n'}(t) - \frac{S}{4n}| \leq \frac{\epsilon}{2n} \]
%for each other $i$ and $i'$,
%\[ |z_i(t) - \frac{S}{2n}| \leq \frac{\epsilon}{2n} \mbox{    and    }  |z_{i'}(t) - \frac{S}{2n}| \leq \frac{\epsilon}{2n}. \]
\[ |z_i(t) - \pi_i F(1)| \leq \frac{\epsilon}{ 2 n}. \]

\noindent {\bf Proof:} \ao{We show that after $O(n \log (n/\epsilon'))$ iterations, the inequality
\begin{equation} \label{Fineq}  |z_i(t) - \pi_i F(1)| \leq \frac{\epsilon'}{2n} F(1) \end{equation} is satisfied. Taking $\epsilon' = \epsilon/F(1)$ in this statement 
yields the lemma.} 

\ao{To prove Eq. (\ref{Fineq}), observe first that scaling the density $\rho$ multiplies both sides by
the same number so that we may assume without loss of generality that $F(1)=1$.}  In this case, $z_i(t)$ is the probability that the random walk that starts at node $k$ with probability $z_k(0)$
is at node $i$ at time $t$. This is a convex combination of the entries of the $i$'th column of $K^t$. By the previous lemma, each entry of the $i$'th entry of each row is not 
more than $\ao{ \epsilon'/(2   n})$ from $\pi_i$ after $O(n \log (n/\epsilon'))$ steps, and the convex combination of these entries must have the same property. {\bf q.e.d.}

\bigskip

\noindent {\bf Proof of Theorem 6.} \aor{Since we have already shown that the positions of the agents have limits which achieve optimal coverage, all that needs to be argued is that every agent is within $\epsilon$ of its final position after $O(n \log (\rho_{\rm max} n/( \rho_{\rm min} \epsilon)))$ rounds}. 

\ao{Since the physical locations of the agents throughout the execution of the dynamic control law do not change if we scale the density $\rho$, we may assume that our density is $\rho'(z) = \rho(z)/\rho_{\rm min}$. Note that $F'(1) = \int_0^1 \rho'(z) dz$ satisfies $F'(1) \leq \rho_{\rm max}/\rho_{\rm min}$. We now apply the previous lemma to the density $\rho'$, and obtain that after   $O(n \log (\rho_{\rm max} n/(\rho_{\rm min} \epsilon)))$ updates, we will
always have} 
\begin{eqnarray*} |z_1(t) + z_{1'}(t) - \frac{F(1)}{2 \aoc{n}} |  & \leq & \frac{  \epsilon}{ n} \\
|z_i(t) + z_{i'}(t) - \frac{F(1)}{\aoc{n}}| & \leq & \frac{ \epsilon}{ n}, ~~~ i = 2,\ldots, n-1 \\ 
|z_n(t) + z_{n'}(t) - \frac{F(1)}{2 \aoc{n}} |  & \leq & \frac{ \epsilon}{ n} 
\end{eqnarray*} which combined with Eq. (\ref{ilimit}) and Eq. (\ref{1limit}) implies that each $z_i(t)+ z_{i'}(t)$ is always within $\epsilon/n$ of its final position after $t_1=\ao{O(n \log (\rho_{\rm max} n/( \rho_{\rm min} \epsilon)))}$ iterations. 

\aoc{Denote the limiting position of agent $i$ by $\overline{x_i}$ and the limiting value of $z_i(t)+z_{i'}(t)$ by $\overline{z_i}$.} \aor{ At some time $t_1' \in [t_1,t_1+U]$ agent $1$ is going to initiatate motion. Then at $t = t_1'+1$,}  \[ |\int_{x_{1}(t)}^{\overline{x_1}} \rho'(u) dz | = |\int_0^{x_1(t)} \rho'(u) du - \int_0^{\overline{x_1}} \rho'(u) du | = |z_1(t-1) + z_{1'}(t-1) - \overline{z_1}| \leq \frac{\epsilon}{n}, \] and since $\rho'(z) \geq 1$, this implies in particular that $|x_1(t) - \overline{x_1}| \leq \epsilon/n$. \aor{Moreover, this remains true for all $t$ after $t_1'+1$; indeed, it does not change when agent $1$ is not initiating motion, and when it is by exactly the above chain of inequalities. }

\aor{We next argue that for all $t \geq t_1+i$}, we have \[ | \int_{x_{i}(t)}^{\overline{x_i}} \rho'(u) du|  \leq \frac{i \epsilon}{n}, \] which similarly implies $|x_i(t) - \overline{x_i}| \leq i \epsilon/n$. We proceed by induction; indeed,  we have already proven the base case of $i=1$ and now assuming we have proven this for $i=l$, we have that at time $t_1'+l+1$,
%\[ |x_{l+1}(t) - \overline{x_{l}}| \leq  |x_{l+1}(t) - x_l(t-1)| + |x_l(t-1) - \overline{x_l}| \leq z_{l+1}(t) + \frac{l \epsilon}{n} %\leq \frac{(l+1) \epsilon}{n}.\] 
\begin{eqnarray*} |\int_{x_{l+1}(t)}^{\overline{x_{l+1}}} \rho'(u) du | 
& = & |\int_{{x_{l+1}(t) }}^{x_{l}(t)} \rho'(u) du  + \int_{x_{l}(t) }^{\overline{x_{l}} } \rho'(u) du  +  \int_{\overline{x_{l}}}^{\overline{x_{l+1}}} \rho'(u) du | \\
& \leq &  |\int_{ x_{l}(t)}^{ \overline{x_{l}}} \rho'(u) du | + |- \int_{{x_{l}(t)}}^{x_{l+1}(t)} \rho'(u) du + \int_{\overline{x_{l}}}^{\overline{x_{l+1}}} \rho'(u) du | \\
& \leq & \frac{l \epsilon}{n} + |- (z_i(t) + z_{i'}(t)) + \overline{z_i}| \\
& \leq & \frac{(l+1) \epsilon}{n}.
\end{eqnarray*} \aor{Moreover, this remains true for all $t \geq t_1' + i$: the inequality is unchanged when agent $i$ is not initiating motions, and holds by 
exactly the same sequence of inequalities when it is. }

Thus after \aor{$t_1+2U$} steps,  every agent is always at most $\epsilon$ away from its limiting location. \aor{Since $U = O(n)$,} this proves the theorem.   {\bf q.e.d.}

\smallskip

\noindent {\bf Remark:} We remark that it is possible to describe a modification of the dynamic control law that is robust to unpredictable \aoj{deletion and addition of agents. We \aoc{informally} sketch this modification here; it is based on the observation that the dynamic control law converges to the correct solution regardless of initial conditions as long as $\sum_{i=1}^n z_i(0) + z_{i'}(0)$ has the correct value of $F(1)$. Consequently, it suffices to describe a protocol which preserves the sum $\sum_{i=1}^n z_i(t) + z_{i'}(t)$ in the event of agent addition and deletion.} 

\aoj{It is easy to deal with the addition of an agent: the new agent simply sets its variables $z_i,z_{i'}$ to zero, thus preserving the sum. On the other hand, if agent $i$ is deleted we must rely on the measurements of neighboring agents to recover $z_i+z_{i'}$. However, observe that the rules of our control law imply that there is a relation between the sum $z_{i}+z_{i'}$ and the \aoc{location of agent $i$. Indeed, agent $i$'s location at the end of update round $t$ is $z_i(t) + z_{i'}(t)$ to the right of the  location of agent $i-1$. Thus the neighbors of agent $i$ can reconstruct $z_i(t) + z_{i'}(t)$ after the unpredictable deletion of agent $i$.  Once this is completed}, the neighboring agents can increase their own variables in any way that preserves the sum $\sum_{i=1}^n z_i(t) + z_{i'}(t)$. }

\aoj{A similar technique may be used to deal with changes in density; if the agents can detect these changes and increase or decrease their own variables $z_i,z_{i'}$ in response, then as long as the number of agents and the density eventually stablize, a suitable modification of the dynamic control law will converge to the correct answer. An interesting open problem is to prove a dynamic guarantee on performance if the number of agents or the density is continually changing. }

\section{SIMULATIONS} We report here on several simulations of our coverage control laws. We are able to observe that quite often the performance is considerably better than the theoretical upper bounds derived in this paper, and
that the dynamic control law of Section IV gives a considerable practical speedup over the static control law of Section
III.  

\begin{figure}[h]
\begin{center}$
\begin{array}{c}
\includegraphics[width=0.51\textwidth]{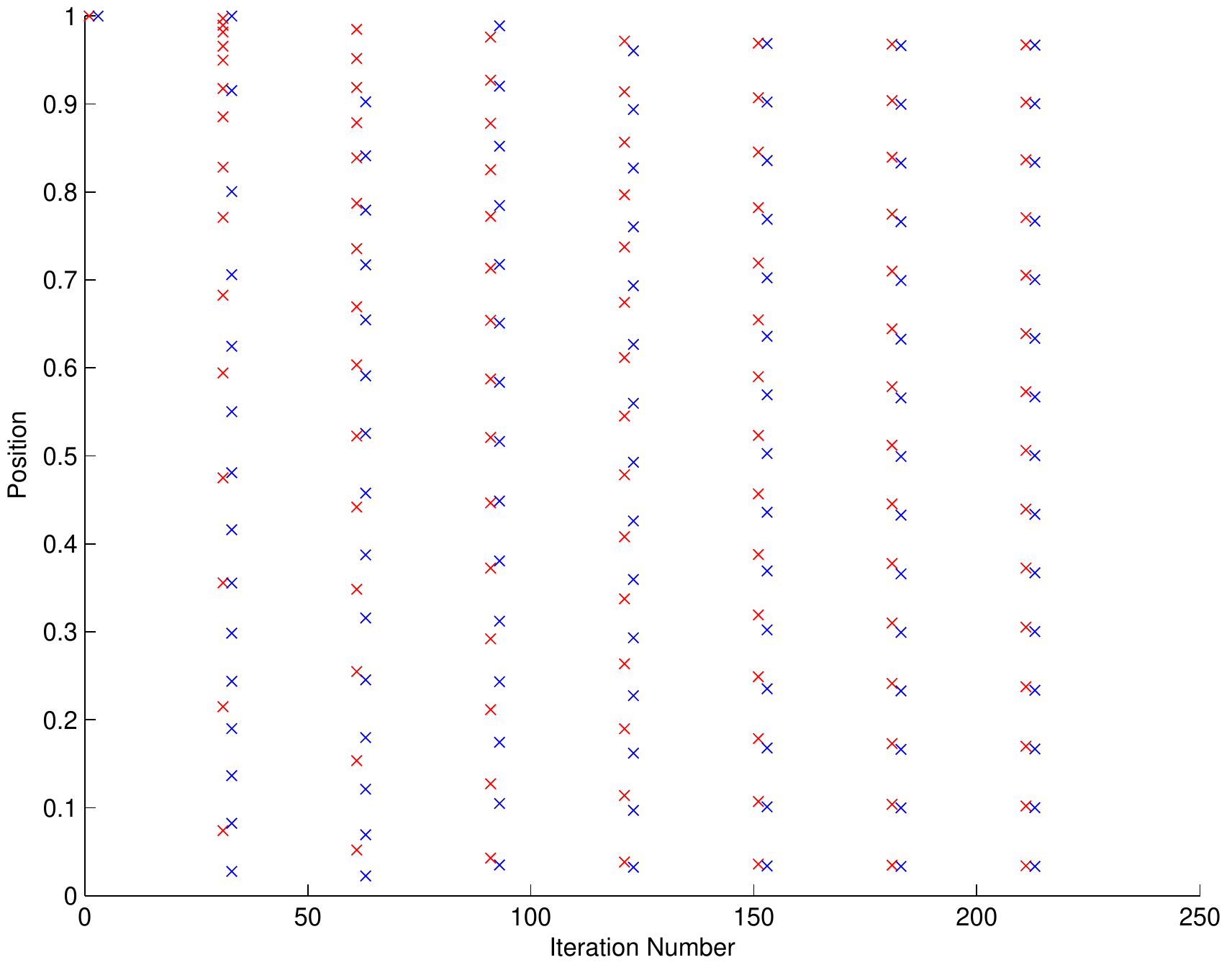} 
\includegraphics[width=0.51\textwidth]{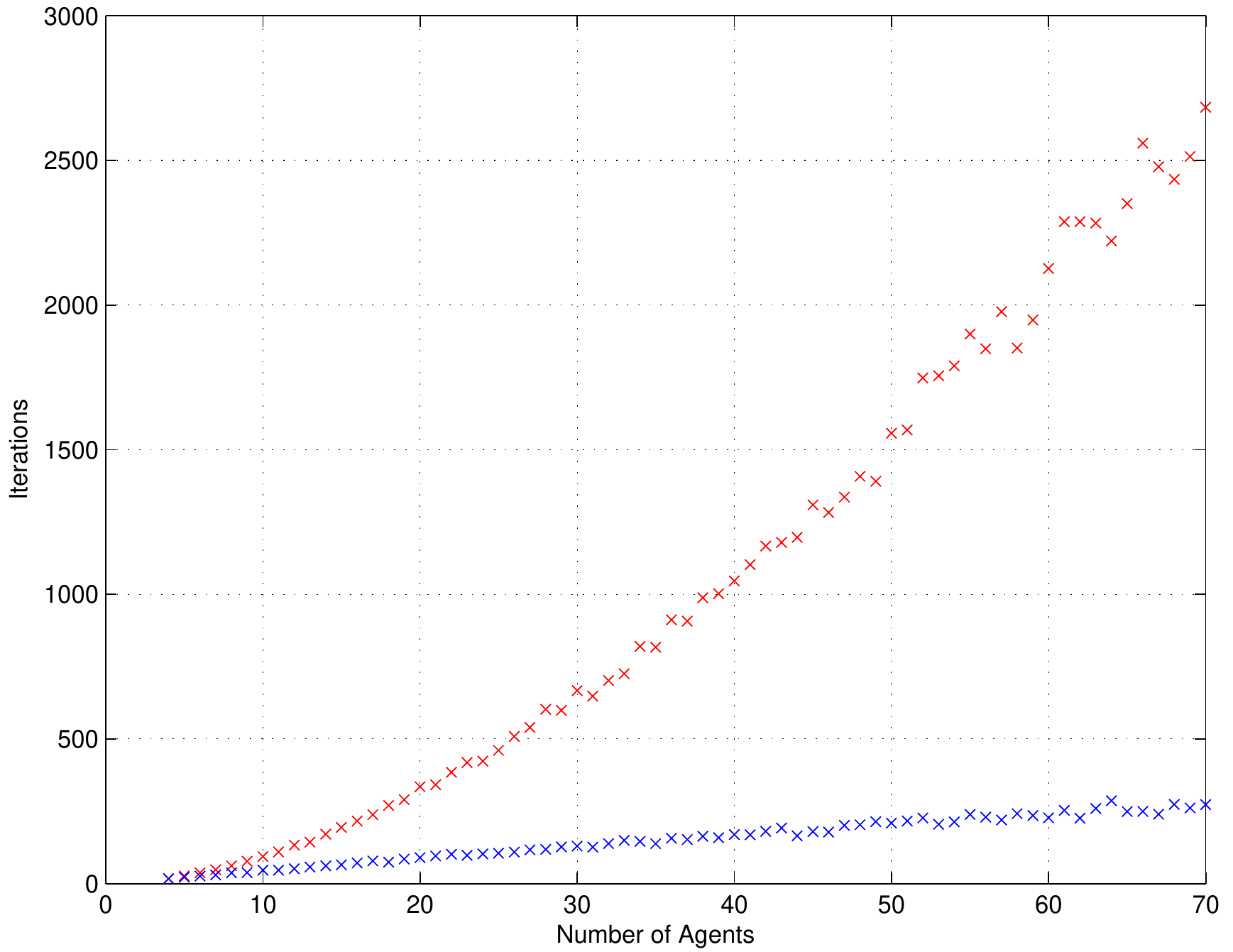}
\end{array}$
\end{center}
\caption{\aoj{Left} figure displays the progress of both the static (red) and dynamic (blue) coverage control laws in the uniform density $\rho(u)=1$. In both cases, the initial condition consists of $n=15$ nodes with the same randomly generated positions. The iteration number $t$ is represented on the x-axis, and the distribution of agents in $[0,1]$ on the y-axis. The \aoj{right figure} shows the number of agents on the x-axis, and the number of iterations until agents are close to their final values on the y-axis. Initial conditions are random as well, and each point on
the right graph is an average of $40$ runs.}
\end{figure}

Figure 3 shows the results from a simulation with random initial conditions. In this case, $x_i(0)$ is the $i$'th largest value of  $n$ random variables all uniform on $[0,1]$. The density $\rho$ is uniform on $[0,1]$. Moreover, we assumed that each agent knows the total number of agents in the system, i.e. $U=n$.  The left figure shows some snapshots from the
progress of both algorithms, while the right figure shows the time until the stopping condition $\sum_{i=1}^n (x_i(t) - \lim_t x_i(t))^2 \leq 10^{-4}$ holds for the first time. The results show reasonably quick convergence for the random initial conditions.  We see that \ao{for the systems shown in the figure}, the static control law has a convergence time which grows slower 
than the quadratic growth proved in Theorem 1, while the dynamic control law has convergence time which appears to grow somewhat 
slower than the linear upper bound of Theorem 6. 

\begin{figure}[h]
\begin{center}$
\begin{array}{c}
\includegraphics[width=0.51\textwidth]{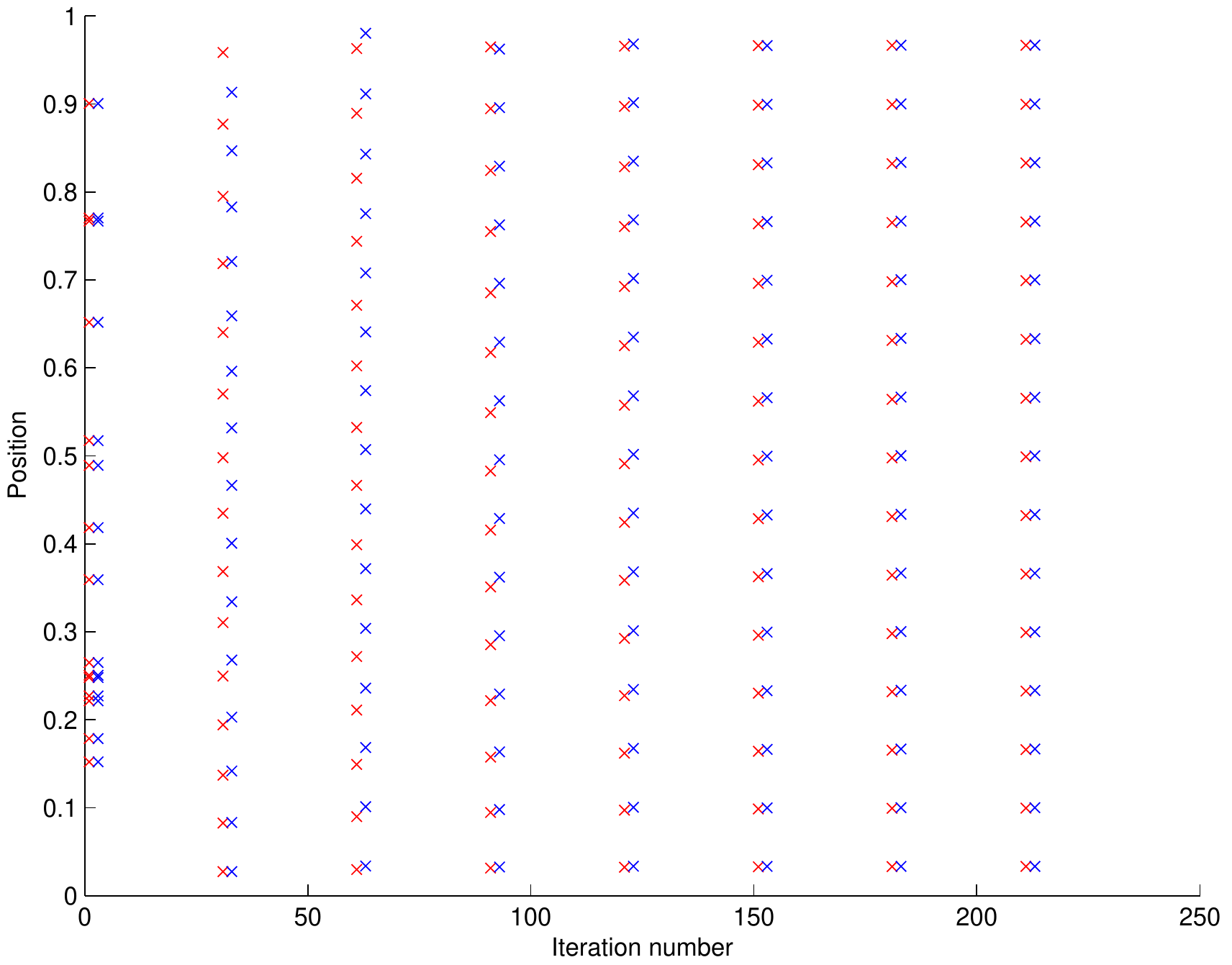} 
\includegraphics[width=0.5\textwidth]{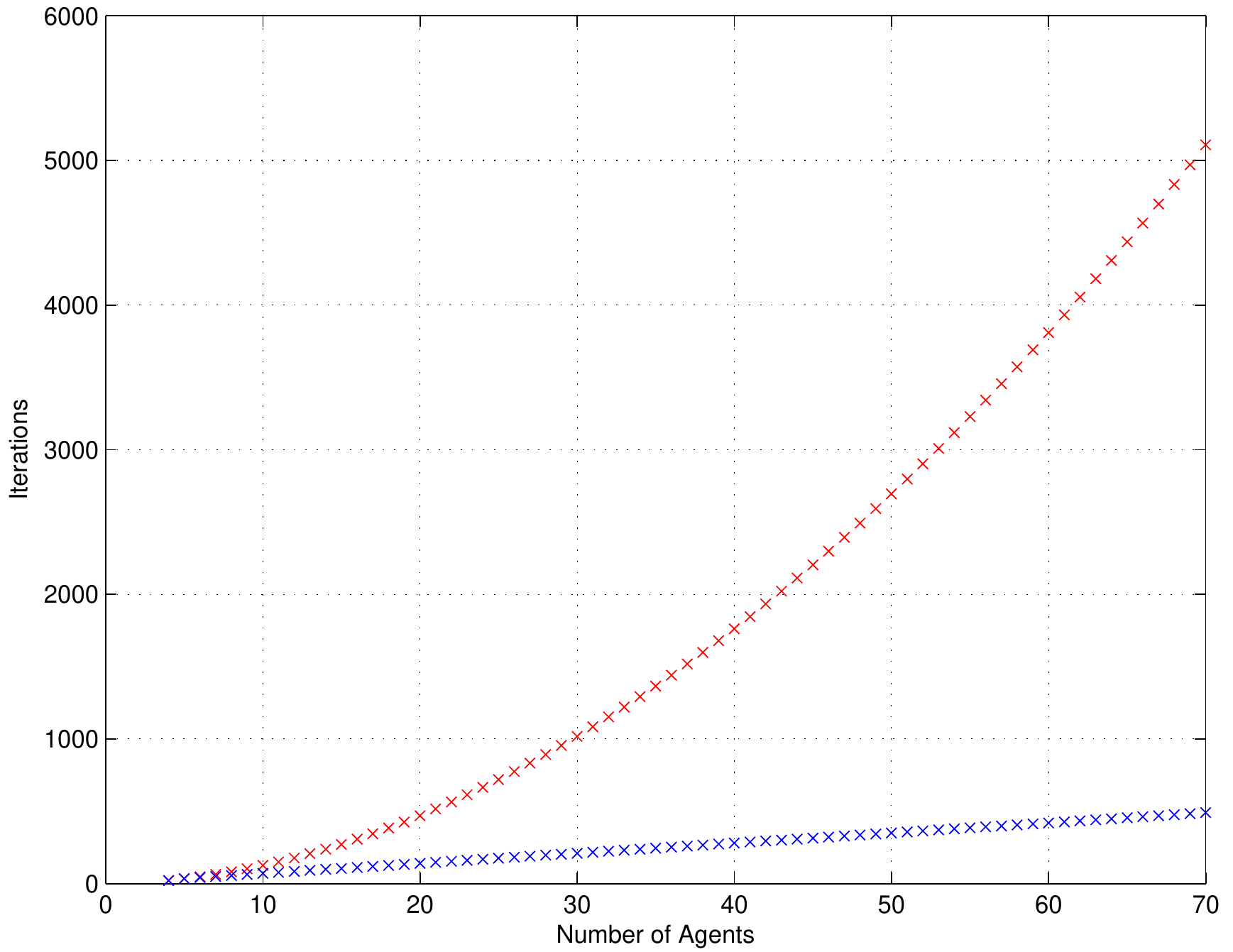}
\end{array}$
\end{center}
\caption{\aoj{Left} figure displays the progress of both the static (red) and dynamic (blue) coverage control laws in the uniform density field $\rho(u)=1$. In both cases, the initial conditions consist of $n=15$ nodes all of which start with $x_i(0)=1$. The iteration number $t$ is represented on the x-axis, and the distribution of agents in $[0,1]$  on the y-axis. The \aoj{right} figure shows the number of agents on the x-axis, and the number of iterations until agents are close to their final values on the y-axis; just as in the left figure, every node starts out with $x_i(0)=1$.}
\end{figure}

On the other hand, \ao{it is not hard to find examples for} which the bounds on growth rates from Theorems 1 and 6 do occur. In Figure 4, we see the performance when every agent starts with $x_i(0)=1$; every other aspect is the same as the simulation in Figure 3. We see that in this case the convergence times do seem to grow quadratically in $n$ for the static control law, and linearly in $n$ for the dynamic control law. 

Intuitively, the issue of convergence rate seems
to be related to how fast information propagates through the network of nodes. The high weights used in the dynamic
control law  - specifically, the weights $1-1/U$, which are quite close to $1$ - allow for fast diffusion of information if an agent is  far from where it
should be. On the other hand, in the static control law information diffuses through nearest neighbor interactions, and this process can be an order of magnitude slower. 

\begin{figure}[h]
\begin{center}$
\begin{array}{c}
\includegraphics[width=0.5\textwidth]{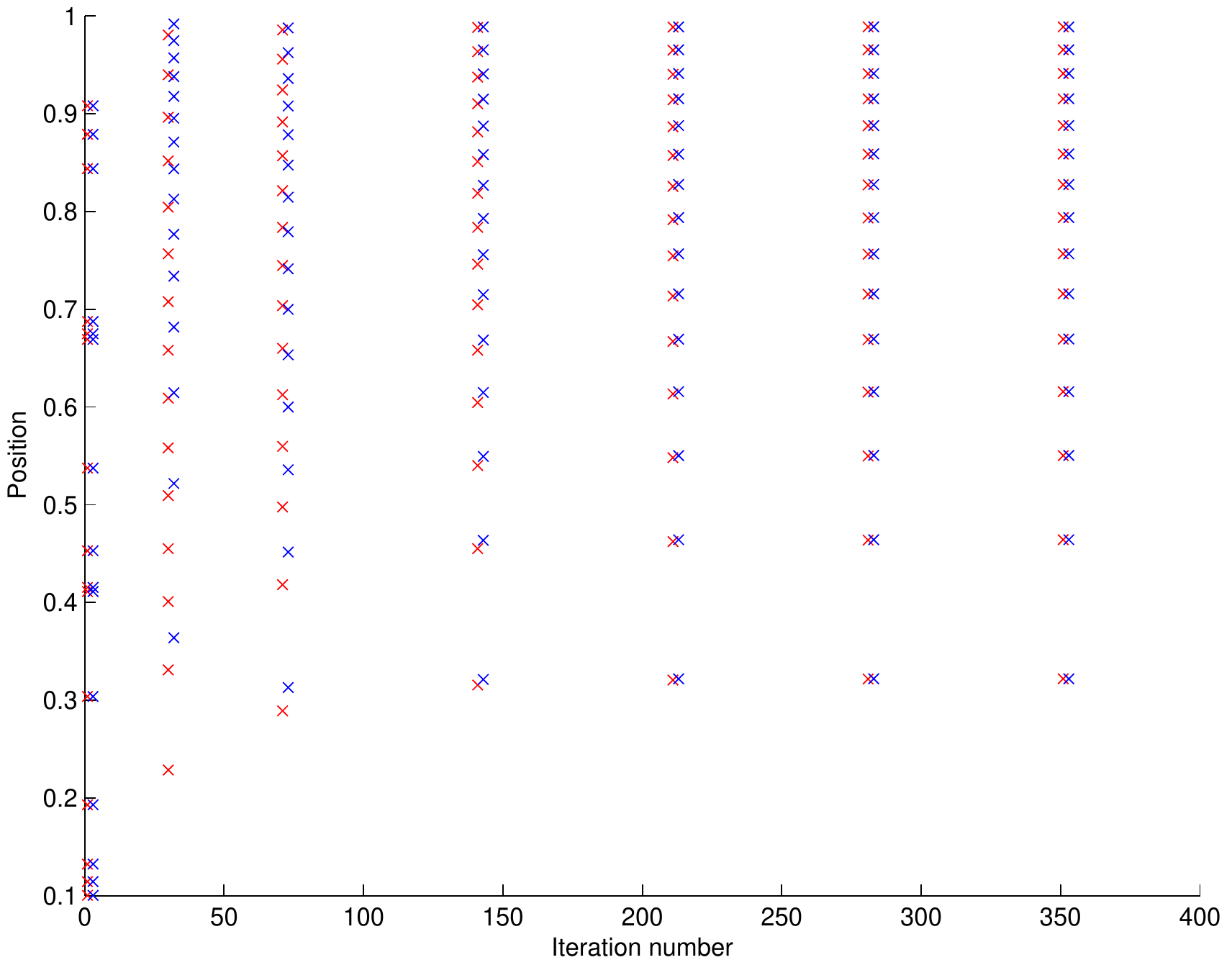} 
\includegraphics[width=0.5\textwidth]{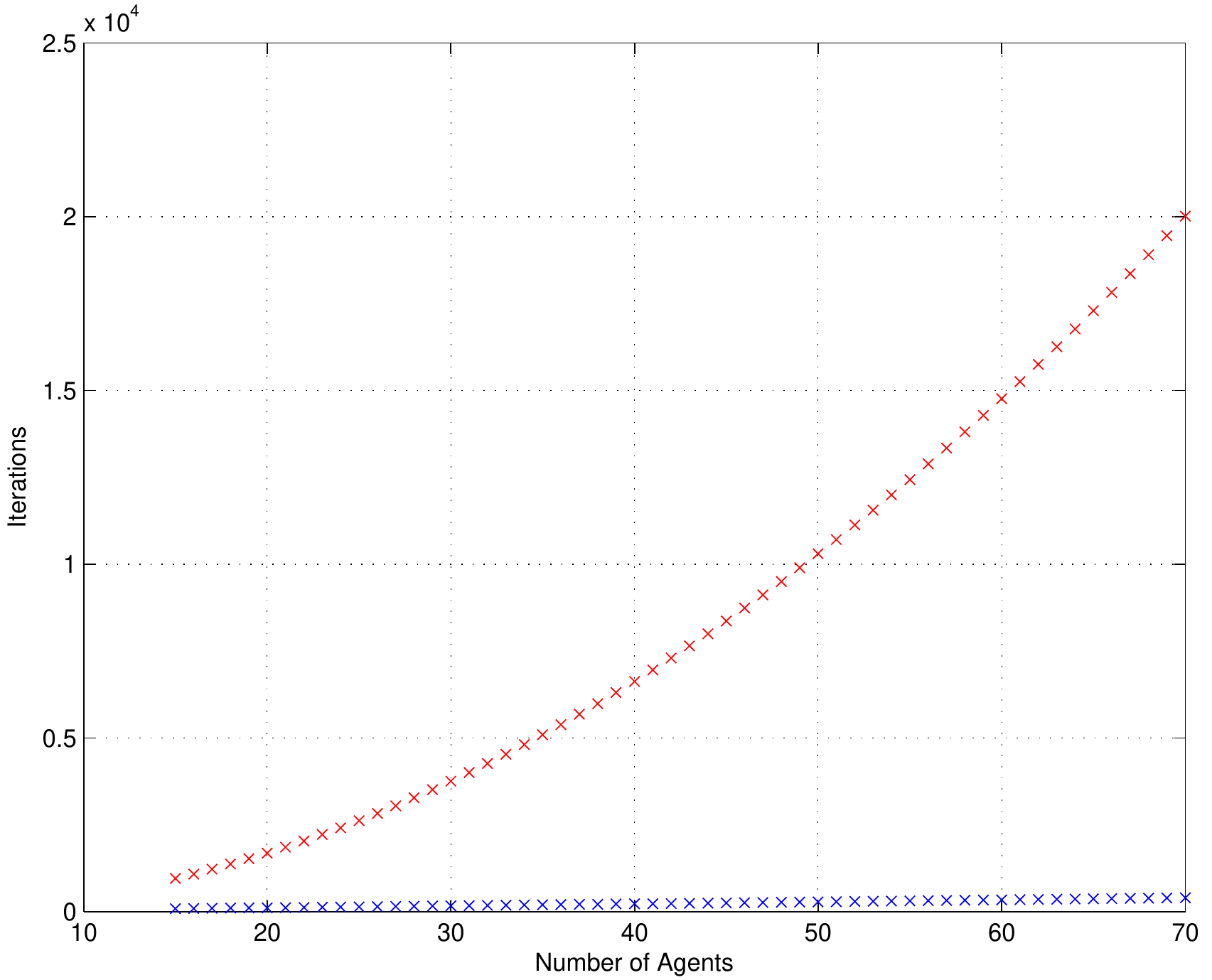}
\end{array}$
\end{center}
\caption{\aoj{The left figure displays the progress of both the static (red) and dynamic (blue) coverage control laws in the nonuniform density field $\rho(x)=x^2$. The initial conditions consist of $n=15$ nodes all of which start with $x_i(0)=0$. The iteration number $t$ is represented on the x-axis, and the distribution of agents in $[0,1]$ on the y-axis. The right figure shows the number of agents on the x-axis, and the number of iterations until agents are close to their final values on the y-axis; just as in the left figure, every node starts out with $x_i(0)=0$.}}
\end{figure}

\aoj{Finally, we describe an example in which the agents achieve optimal coverage with respect to a nonuniform metric.  Figure 5 shows a simulation with density $\rho(x)=x^2$, which requires the nodes to concentrate themselves closer to the boundary of $1$ of the interval $[0,1]$. The starting point was chosen to be $x_i(0)=0$ for $i=1,\ldots,n$, in part because it is a difficult starting point for the coverage algorithms we have presented here: the agents will have to move right one-by-one.  Convergence times seem to be somewhat larger here than in the uniform density case, and the gap between the convergence times of the dynamic and static algorithms seems to be more pronounced. The convergence times shown in the graphs are consistent with the bounds of Theorems 1 and 6. The dynamic control law exhibits some interesting  dynamics, first moving many agents too close before pushing them apart to the optimal configuration. It is of interest to understand the non-asymptotic dynamics of this control law better. }

\section{CONCLUSIONS} We have investigated distributed control laws for mobile, autonomous agents to position themselves on the line for 
optimal coverage in a nonuniform field. Our main results are stated in Theorems 1 and 6. Theorem 1 gives a quantitative upper 
bound on the convergence time of a relatively simple control law for coverage. Theorem 6 discusses a dynamic control law
which, while making stronger assumptions on the capabilities of each agent, manages to accomplish the coverage task \aoj{with an order of magntude 
fewer sensing/communication/movement iterations in the worst case.} 

Our work suggests a number of open questions. It is of interest to understand whether the increased capabilities of the agents in Section 4 are really necessary to achieve better performance. In addition, it would be interesting to explore whether the results described here extend to two and higher dimensions, and in particular, whether a dynamic control law such as the one in Theorem 6 might be useful for speeding up performance in more general settings. \aor{Generalizing these results to higher dimensions is challenging because the notions of a ``leftmost'' and ``rightmost'' agent does not
make sense in two dimensions; moreover, our algorithm crucially relies on the transformation $x \rightarrow \int_0^x \rho(z) ~dz$ which ``stretches'' the domain to make it  uniform; it is not immediately clear how to compute the corresponding transformation in two dimensions in a distributed way.} %\aoj{It would also %be interesting to describe the ability of collections of agents to achieve good coverage when the density is constantly time-varying.} 
\aoc{Finally, an open problem is to replace the assumption that nodes are capable of computing integrals of the density $\rho$ with the weaker assumption that noisy samples of $\rho$ are available to the nodes. }

\section{ACKNOWLEDGEMENTS} 

The authors would like to thank Hari Narayanan for useful discussions.

%%%%%%%%%%%%%%%%%%%%%%%%%%%%%%%%%%%%%%%%%%%%%%%%%%%%%%%%%%%%%%%%%%%%%%%%%%%%%%%%

\end{document}